\documentclass[a4paper]{article}
\usepackage[T1]{fontenc}
\usepackage{amsfonts,amsmath,amsthm,amssymb}
\usepackage{calc}
\usepackage{gensymb}
\usepackage[nottoc, notlof, notlot]{tocbibind}
\usepackage{fullpage}
\usepackage{mathrsfs}  
\usepackage{bigints}
\usepackage{indentfirst}
\usepackage{enumerate}
\usepackage{authblk}
\usepackage{hyperref}

\DeclareMathOperator{\tr}{Tr}
\DeclareMathOperator{\op}{op}
\DeclareMathOperator{\cov}{Cov}
\DeclareMathOperator{\var}{Var}

\DeclareMathOperator{\supp}{Supp}
\let\limsup\relax
\DeclareMathOperator*{\limsup}{limsup}
\let\liminf\relax
\DeclareMathOperator*{\liminf}{liminf}

\newcommand{\norm}[1]{\left\Vert #1\right\Vert}

\begin{document}

\newtheorem{theorem}{Theorem} [section]
\newtheorem{prop}{Proposition} [section]
\newtheorem{defi}{Definition} [section]
\newtheorem{exe}{Exemple} [section]
\newtheorem{lemma}{Lemma} [section]
\newtheorem{rem}{Remark} [section]
\newtheorem{cor}{Corollary} [section]
\newtheorem{conj}{Conjecture} [section]
\renewcommand\P{\mathbb{P}}
\newcommand\E{\mathbb{E}}
\renewcommand\u{\mathbf{u}}
\renewcommand\v{\mathbf{v}}
\newcommand\N{\mathbb{N}}
\newcommand\1{\mathbf{1}}
\newcommand\C{\mathbb{C}}
\newcommand\CC{\mathcal{C}}
\newcommand\M{\mathbb{M}}
\newcommand\R{\mathbb{R}}
\newcommand\U{\mathbb{U}}
\newcommand\A{\mathcal{A}}
\newcommand\PP{\mathcal{P}}
\newcommand\B{\mathcal{B}}
\newcommand\D{\mathcal{D}}
\renewcommand\i{\mathbf{i}}
\renewcommand\S{\mathcal{S}_{N,t}}
\renewcommand\L{\mathcal{L}_{\Phi}}
\renewcommand\d{\partial_i}
\renewcommand\.{\ .}
\renewcommand\,{\ ,}

\def\etc{,\dots ,}

\title{On the operator norm of non-commutative polynomials in deterministic matrices and iid Haar unitary matrices}

\date{}

\author[1]{F\'elix Parraud}
\affil[1]{\small Universit\'e de Lyon, ENSL, UMPA, 46 all\'ee d'Italie, 69007 Lyon. Department of Mathematics, Graduate School of Science, Kyoto University, Kyoto 606-8502, Japan.}

\maketitle
\noindent E-mail of the corresponding author: \href{mailto:felix.parraud@ens-lyon.fr}{felix.parraud@ens-lyon.fr}

\noindent Data sharing not applicable to this article as no datasets were generated or analysed during the current study.

\begin{abstract}
	
	Let $U^N = (U_1^N\etc U^N_p)$ be a $p$-tuple of $N\times N$ independent Haar unitary matrices and $Z^{NM}$ be any family of 
	deterministic matrices in $\M_N(\C)\otimes \M_M(\C)$. Let $P$ be a self-adjoint non-commutative polynomial. In \cite{voichaar}, Voiculescu showed that the empirical measure of the eigenvalues of this polynomial evaluated in Haar unitary matrices and deterministic matrices converges towards a deterministic measure defined thanks to  free probability theory. Now, let $f$ be a smooth function. The main technical result of this paper is a precise bound of the difference between the expectation of
	$$ \frac{1}{MN} \tr_{\M_N(\C)}\otimes\tr_{\M_M(\C)}\left( f(P(U^N\otimes I_M,Z^{NM})) \right) , $$
	
	\noindent and its limit when $N$ goes to infinity. If $f$ is seven times differentiable, we show that it is bounded by $M^2 \norm{f}_{\mathcal{C}^6} \ln^2(N)\times N^{-2}$. As a corollary we obtain a new proof with quantitative bounds of a result of Collins and Male which gives sufficient conditions for the operator norm of a polynomial evaluated in Haar unitary matrices and deterministic matrices to converge almost surely towards its free limit. Our result also  holds in much greater generality. For instance,  it allows to prove that  if $U^N$ and $Y^{M_N}$ are independent and $M_N=o(N^{1/3}\ln^{-2/3}(N))$, then the norm of any polynomial in $(U^N\otimes I_{M_N}, I_N\otimes Y^{M_N})$ converges almost surely towards its free limit. Previous results required that $M=M_N$ is constant.
	
\end{abstract}

\section{Introduction}

Understanding the behaviour of random matrices in large dimension is the core of random matrix theory. In the early nineties Voiculescu showed that one could get very accurate results with the help of non-commutative probability theory. This theory is equipped with a notion of freeness, analogous to independence in classical probability theory, which often describes accurately the asymptotic behaviour of random matrices. In \cite{Vo91} he studied the asymptotic behaviour of independent  matrices taken from the Gaussian Unitary Ensemble (GUE). In a later paper he proved a similar theorem for Haar unitary matrices, which are random matrices whose law is the Haar measure on the unitary group $\U_N$. In a nutshell, Voiculescu proved in \cite{voichaar} that given $U^N_1,\dots,U^N_p$ independent Haar unitary matrices, the renormalized trace of a polynomial $P$ evaluated in these matrices converges towards a deterministic limit $\alpha(P)$. 
Specifically, the following holds true almost surely:

\begin{equation}\label{2dv} \lim_{N\to \infty} \frac{1}{N}\tr_N\left( P(U_1^N,\dots,U_p^N,{U_1^N}^*,\dots,{U_p^N}^*) \right) = \alpha(P) .\end{equation}

\noindent Voiculescu computed the limit $\alpha(P)$ with the help of free probability. To give more detail, let $B_N$ be a self-adjoint matrix of size $N$, then one can define the empirical measure of its (real) eigenvalues by 

$$ \mu_{B_N} = \frac{1}{N} \sum_{ i=1}^N \delta_{\lambda_i} ,$$

\noindent where $\delta_{\lambda}$ is the Dirac mass in $\lambda$ and $\lambda_1\etc \lambda_N$ are the eingenvalue of $B_N$. In particular, if $P$ is a self-adjoint polynomial, that is such that for any matrices $A_1\etc A_d$, $P(A_1\etc A_d,A_1^*\etc A_d^*)$ is a self-adjoint matrix, then one can define the random measure $\mu_{P(U_1^N,\dots,U_p^N,{U_1^N}^*,\dots,{U_p^N}^*)}$. In this case, Voiculescu's result \eqref{2dv} implies that there exists a measure $\mu_P$ with compact support such that almost surely $\mu_{P(U_1^N,\dots,U_p^N,{U_1^N}^*,\dots,{U_p^N}^*)}$ converges weakly towards $\mu_P$: its moments are given by $\mu_P(x^k)=\alpha(P^k)$ for all integer numbers $k$.\\

However, the convergence of the empirical measure of the eigenvalues of a matrix does not say much about the local properties of its spectrum, in particular about the convergence of the norm of this matrix, or the local fluctuations of its spectrum. For a comprehensive survey of important milestones related to these questions, we refer to the introduction of our previous paper \cite{un}. In a nutshell, when dealing with a single matrix, incredibly precise results are known.  Typically, concerning the GUE, very precise results were obtained by Tracy and Widom in the early nineties in \cite{TW1}. On the other hand, there are much less results available when one deals with a polynomial in several random matrices. One of the most notable result was found by  Haagerup and Thorbj\o rnsen in 2005 in \cite{HT}: they proved the almost sure convergence of the norm of a polynomial evaluated in independent GUE matrices. Equivalently, for $P$ a self-adjoint polynomial, they proved that almost surely, for $N$ large enough, 
\begin{equation}
\label{2spec}
\sigma\left( P(X_1^N,\dots,X_d^N) \right) \subset \supp \mu_P + (-\varepsilon,\varepsilon) ,
\end{equation}

\noindent where $\sigma(H)$ is the spectrum of $H$ and $\supp \mu_P$ the support of the measure $\mu_P$. The result \eqref{2spec} was a major breakthrough in the context of free probability and was refined in multiple ways, see \cite{SCH,CD,AND,belin-capi,pisier,un}. Those results all have in common that the basic random matrix is always self-adjoint. Much less is known in the non self-adjoint case. However Collins and Male proved in \cite{collins_male} the same result as in \cite{male} but with unitary Haar matrices instead of GUE matrices by using Male's former paper. With the exception of \cite{collins_male} and \cite{un}, all of these results are essentially based on the method introduced by Haagerup and Thorbj\o rnsen who relies on the so-called linearization trick. The main idea of this tool is that given a polynomial $P$,  the spectrum of $P(X_1^N,\dots,X_d^N)$ is closely related to  the spectrum of 
$$ L_N = a_0\otimes I_N + \sum_{ i=1}^d a_i\otimes X^N_i ,$$

\noindent where $a_0\etc a_d$ are matrices of size $k$ depending only on $P$. Thus we trade a polynomial of degree $d$ with coefficient in $\C$ by a polynomial of degree $1$ with coefficient in $\M_{k(d)}(\C)$. In \cite{collins_male}, the main idea was to view Haar unitary matrices as a random function of a GUE random matrix. Then the authors showed that almost surely this random function converges uniformly and they concluded by using the main result of \cite{male}. An issue of this method is that it does not give any quantitative estimate. An important aim of this paper is to remedy to this problem. Our approach requires neither the linearization trick, nor the study of the Stieljes transform and attacks the problem directly without using previous results about the strong convergence of GUE random matrices. In this sense the proof is more direct and less algebraic. We will apply it to a generalization of Haar unitary matrices by tackling the case of Haar unitary matrices tensorized with deterministic matrices.

A usual strategy to study outliers, that are the eigenvalues going away from the spectrum,  is to study the \emph{non-renormalized} trace of smooth non-polynomial functions evaluated in independent Haar matrices i.e. if $P$ is self-adjoint:
$$\tr_N\left( f\left(P\left(U_1^N,\dots,U_p^N,{U_1^N}^*,\dots,{U_p^N}^*\right)\right) \right) .$$
This strategy was also  used by Haagerup, Thorbj\o rnsen and Male. Indeed it is easy to see that if $f$ is a function which takes value $0$ on $(-\infty,C-\varepsilon]$, $1$ on $[C,\infty)$ and in $[0,1]$ 
elsewhere, then with $\lambda_1(X)$ the largest eigenvalue of $X$,
$$ \P\Big( \lambda_1(P(U_1^N,\dots,U_p^N,{U_1^N}^*,\dots,{U_p^N}^*)) \geq C \Big) \leq\ \P\Big( \tr_{N}\left( f(P(U_1^N,\dots,U_p^N,{U_1^N}^*,\dots,{U_p^N}^*)) \right) \geq 1 \Big) .$$

\noindent Hence, if we can prove that $\tr_{N}\left( f(P(U_1^N,\dots,U_p^N,{U_1^N}^*,\dots,{U_p^N}^*)) \right)$ converges towards $0$ in probability, this would already 
yield expected results. The above is just a well-known example, but one can get much more out of this strategy. 
Therefore, we need to study the non-renormalized trace. The case where $f$ is a polynomial function has already been studied a long time ago, starting with the pioneering works \cite{CDu,GTCL}, and later formalized by the concept of second order freeness \cite{MingoSpeicher,MingoSpeicher2}. However here we have to deal with a function $f$ which is at best $C^{\infty}$. This makes things considerably more difficult and forces us to adopt a completely different approach. The main result is the following theorem (for the notations, we refer to Section \ref{2definit} -- for now, let us specify that $\frac{1}{N}\tr_N$ denotes the usual renormalized trace on $N\times N$ matrices whereas $\tau$ denotes its free limit):

\begin{theorem}
	\label{2imp1}
	We define
	\begin{itemize}
		\item $u=(u_1,\dots ,u_p,u_1^*,\dots ,u_p^*)$ a family of $p$ free Haar unitaries and their adjoints,
		\item $U^N = (U_1^N,\dots,U_p^N,(U_1^N)^*,\dots,(U_p^N)^*)$ i.i.d. Haar unitary matrices of size $N$, and their adjoints.
		\item $Z^{NM} = (Z_1^{NM},\dots,Z_q^{NM})\in\M_N(\C)\otimes \M_M(\C)$ deterministic matrices and their adjoint,
		\item $P$ a self-adjoint polynomial,
		\item $f:\R\mapsto\R$ a function at least six times differentiable.
	\end{itemize}
	
	\noindent Then there exists a polynomial $L_P\in\R^+[X]$ which only depends on $P$ such that for any $N,M$,
	
	\begin{align*}
	\Bigg| &\E\left[\frac{1}{MN}\tr_{MN}\Big(f\left(P\left(U^N\otimes I_M,Z^{NM}\right)\right)\Big)\right] - \tau_N\otimes\tau_M\Big(f\left(P\left(u\otimes I_M,Z^{NM}\right)\right)\Big) \Bigg| \\
	&\leq \frac{\ln^2(N) M^2}{N^2} L_P\left(\norm{Z^{NM}}\right)\times \norm{f}_{\mathcal{C}^6}.
	\end{align*}
	
	\noindent where $\norm{Z^{NM}} = \sup\limits_{1\leq i\leq q} \norm{Z_i^{NM}}$ and $\norm{f}_{\mathcal{C}^k}$ is the sum of the supremum on $\R$ of the first $k$ derivatives. If $Z^{NM} = (I_N\otimes Y_1^M,\dots ,I_N\otimes Y_q^M)$ and that these matrices commute, then we have the same inequality without the $M^2$.
	
\end{theorem}

This theorem is a consequence of the slightly sharper, but less explicit, Theorem \ref{2imp2}. Those two theorems are essentially the same, but in Theorem \ref{2imp2}, instead of having the norm $C^6$ of $f$, we have the fourth moment of the Fourier transform of $f$. The above theorem calls for a few remarks. 

\begin{itemize}
	\item We assumed that the matrices $Z^{NM}$ were deterministic, but thanks to Fubini's theorem we can assume that they are random 
	matrices as long as they are independent from $U^N$. In this situation though, $L_P\left(\norm{Z^{NM}}\right)$ in the right side of the inequality is a random variable (and thus we need some additional assumptions on the law of $Z^{NM}$ if we want its expectation to be finite for instance).
	\item In Theorems \ref{2imp1} and \ref{2imp2} we have $U^N\otimes I_M$ and $u\otimes I_M$, however it is very easy to replace them by $U^N\otimes Y^M$ and $u\otimes Y^M$ for some matrices $Y^M_i\in M_M(\C)$. Indeed we just need to apply Theorem \ref{2imp1} or 
	\ref{2imp2} with $Z^{NM}=I_N\otimes Y^M$. Besides, in this situation, $L_P\left(\norm{Z^{NM}}\right) =  L_P\left(\norm{Y^M}\right)$ does not depend on $N$. What this means is that if we have a matrix whose coefficients are polynomial in $U^N$, and that we replace $U^N$ by $u$, we only change the spectra of this matrix by $M^2N^{-2}$ in average.
	\item In the specific case where $Z^{NM} = (I_N\otimes Y_1^M,\dots ,I_N\otimes Y_q^M)$ and the $Y_i^M$ commute, as we stated in Theorem \ref{2imp1}, we have the same inequality without the $M^2$. Lowering the exponent in all generality would yield a direct improvement to Theorem \ref{2strongconv} but we currently do not know whether it is actually possible. A lead to do so would be to prove a sharper version of Lemma \ref{2normineq}. While this seems unrealistic for deterministic matrices, it might be possible to get some results when the matrices $Y_i^M$ are random.
\end{itemize}

A detailed overview of the proof is given in Subsection \ref{2overview}. Similarly to \cite{un}, we interpolate Haar unitary matrices and free Haar unitaries with the help of a free Ornstein-Uhlenbeck process on the unitary group, i.e. the free unitary Brownian motion. For a reference, see Definition \ref{2unitarybro}. However in \cite{un} this idea was only to understand the intuition of the proof. In this paper the computations involved were quite different, indeed since we were considering the usual free Ornstein-Uhlenbeck process, we could use a computation trick to replace this process by a well-chosen interpolation between GUE matrices and free semicirculars. This means that we did not need to use free stochastic calculus. There is no such trick for the free unitary Brownian motion, hence the computations use much more advanced tools.

When using this process, the Schwinger-Dyson equations, which can be seen as an integration by part, appear in the computation. For more information about these equations we refer to \cite{aliceSD} to find numerous applications. In the specific case of the unitary group it is worth checking the proof of Theorem 5.4.10 from \cite{alice}. Even though those equations only come into play in the proof of Lemma \ref{2grcalc}, they play a major role in the proof since we could get a theorem similar to Theorem \ref{2imp1} for any random matrices which satisfies those equations. 

Theorem \ref{2imp1} is the crux of the paper and allows us to deduce many corollaries. Firstly we get the following result. The first statement is basically Theorem 1.4 from \cite{collins_male}. The second one is entirely new and let us tensorize by matrices whose size goes to infinity when until now we could only work with tensor of finite size. This theorem is about strong convergence of random matrices, that is the convergence of the norm of polynomials in these matrices, see Definition \ref{2freeprob}.

\begin{theorem}
	\label{2strongconv}
	Let  the  following objects be given:
	\begin{itemize}
		\item $U^N = (U_1^N,\dots,U_p^N)$ independent unitary Haar matrices of size $N$,
		\item $u = (u_1,\dots,u_p)$ a system of free Haar unitaries,
		\item $Y^M = (Y_1^M,\dots,Y_q^M)$ random matrices of size $M$, which almost surely, as $M$ goes to infinity, converges strongly in distribution towards a $q$-tuple $y$ of non-commutative random variables in a $\mathcal C^*$- probability space $\B$ with a faithful trace $\tau_{\B}$.
		\item $Z^N = (Z_1^N,\dots,Z_q^N)$ random matrices of size $N$, which almost surely, as $N$ goes to infinity, converges strongly in distribution towards a $q$-tuple $z$ of non-commutative random variables in a $\mathcal C^*$- probability space  with a faithful trace,
	\end{itemize}
	
	\noindent then the following holds true.
	
	\begin{itemize}
		\item If  $U^N$ and $Z^N$ are independent, almost surely, $(U^N, Z^N)$ converges strongly in distribution towards 
		$\mathcal{F} = (u,z)$, where $\mathcal{F}$ belongs to a $\mathcal C^*$- probability space $(\A,*,\tau_{\A},\norm{.})$ in which $u$ and $z$ are free.
		\item If $(M_N)_{N\geq 0}$ is a sequence of integers such that $M_N = o(N^{1/3} \ln^{-2/3}(N))$, $U^N$ and $Y^{M_N}$ are independent, then  almost surely $(U^N\otimes I_{M_N}, I_N\otimes Y^{M_N})$ converges strongly in distribution towards $\mathcal{F} = (u\otimes 1, 1\otimes y)$ when $N$ goes to infinity. The family $\mathcal{F}$ thus belongs to $\A\otimes_{\min}\B$ (see Definition \ref{2mini}). Besides if the matrices $Y^{M_N}$ commute, then we can weaken the assumption on $M_N$ by only assuming that $M_N = o(N \ln^{-2}(N))$.
	\end{itemize}
	
\end{theorem}

Understanding the Stieljes transform of a matrix gives a lot of information about its spectrum. This was actually a very important point in the proof of Haagerup and Thorbj{\o}rnsen's theorem. Our proof does not use this tool, however our final result, Theorem \ref{2imp2}, allows us to deduce the following estimate. Being given a self-adjoint $NM\times NM$ matrix, we denote by $G_A$ its Stieltjes transform:
$$G_A(z)=\frac{1}{NM}\tr_{NM}\left(\frac{1}{z-A}\right).$$
This definition extends to the tensor product of free Haar unitaries with deterministic matrices by replacing $(NM)^{-1}\tr_{NM}$ by $\tau_N\otimes \tau_M$.

\begin{cor}
	\label{2stieljes}
	
	Given
	\begin{itemize}
		\item $u=(u_1,\dots ,u_p,u_1^*,\dots ,u_p^*)$ a family of $p$ free Haar unitaries and their adjoints,
		\item $U^N = (U_1^N,\dots,U_p^N,(U_1^N)^*,\dots,(U_p^N)^*)$ i.i.d. Haar unitary matrices of size $N$, and their adjoints.
		\item $Y^M = (Y_1^M,\dots,Y_q^M,{Y_1^M}^*,\dots,{Y_q^M}^*)$ deterministic matrices of size $M$ and their 
		adjoints,
		\item $P$ a self-adjoint polynomial,
	\end{itemize}
	\noindent there exists a polynomial $L_P$ such that for every $Y^M$, $z\in \C\backslash \R$, $M,N\in \N$,	
	$$ \left| \E\left[ G_{P(U^N\otimes I_M, I_N\otimes Y^M)}(z) \right] - G_{P(u\otimes I_M, 1\otimes Y^M)}(z) \right| \leq L_P\left(\norm{Y^M}\right) \frac{M^2 \ln^2(N)}{N^2} \left( \frac{1}{ \left| \Im(z)\right|^5} + \frac{1}{ \left| \Im(z)\right|^2} \right) .$$
	
	\noindent where $\norm{Y^{M}} = \sup\limits_{1\leq i\leq q} \norm{Y_i^{M}}$.
\end{cor}

One of the limitation of Theorem \ref{2imp1} is that we need to pick $f$ regular enough. Actually by approximating $f$, we can afford to take $f$ less regular at the cost of a slower speed of convergence. In other words, we trade some degree of regularity on $f$ for a smaller exponent in $N$. The best that we can achieve is to take $f$ Lipschitz. Thus it makes sense to introduce the Lipschitz-bounded metric which is the standard metric to metrize the topology of the weak convergence of probability measures on $\R$. Let $\mathcal{F}_{LU}$ be the set of Lipschitz function from $\R$ to $\R$, uniformly bounded by $1$ and with Lipschitz constant at most $1$, then
$$ d_{LU}(\mu,\nu) = \sup_{f\in \mathcal{F}_{LU}} \left| \int_{\R} f d\mu - \int_{\R} f d\nu \right| .$$

\noindent This metric is a slight weakening of the Wasserstein-1 distance which is defined similarly but without the assumption of boundedness on the functions $f$. For more information about this metric we refer to Appendix C.2 of \cite{alice}. In this paper, we get the following result:

\begin{cor}
	\label{2LU}
	
	Under the same notations as in Corollary \ref{2stieljes}, there exists a polynomial $L_P$ such that for every matrices $Y^M$ and $M,N\in \N$,	
	$$ d_{LU}\left(\E[\mu_{P(U^N\otimes I_M, I_N\otimes Y_M)}] , \mu_{P(u\otimes I_M, 1\otimes Y_M)}\right) \leq L_P\left(\norm{Y^M}\right) M^2 \left( \frac{\ln N}{N} \right)^{1/3} .$$
\end{cor}

This paper is organized as follows. In section \ref{2definit} we give many usual definitions and notations in free probability, commutative and non-commutative stochastic calculus. Section \ref{2preli} contains the proof of many important properties which we will need later on. Section \ref{2mainsec} contains the proof of Theorem \ref{2imp1}. Finally in section \ref{2proofcoro} we prove all of the corollaries.

\section{Framework and standard properties}
\label{2definit}

\subsection{Usual definitions in free probability}
\label{2deffree}

In order to be self-contained, we begin by reminding the following definitions of free probability.

\begin{defi}
	\label{2freeprob}
	\begin{itemize}
		\item A $\mathcal{C}^*$-probability space $(\A,*,\tau,\norm{.}) $ is a unital $\mathcal{C}^*$-algebra $(\A,*,\norm{.})$ endowed with a state $\tau$, i.e. a linear map $\tau : \A \to \C$ satisfying $\tau(1_{\A})=1$ and $\tau(a^*a)\geq 0$ for all $a\in \A$. In this paper we always assume that $\tau$ is a trace, i.e. that it satisfies $\tau(ab) = \tau(ba) $ for any $a,b\in\A$. An element of $\A$ is called a (non commutative) random variable. We will always work with faithful trace, that is such that if $a\in\A$, $\tau(a^*a)=0$ if and only if $a=0$, in which case the norm is determined by $\tau$ thanks to the formula:
		$$ \norm{a} = \lim_{k\to\infty} \big(\tau\big((a^*a)^{2k}\big)\big)^{1/2k}. $$
		
		\item Let $\A_1,\dots,\A_n$ be $*$-subalgebras of $\A$, having the same unit as $\A$. They are said to be free if for all $k$, for all $a_i\in\A_{j_i}$ such that $j_1\neq j_2$, $j_2\neq j_3$, \dots , $j_{k-1}\neq j_k$:
		$$ \tau\Big( (a_1-\tau(a_1))(a_2-\tau(a_2))\dots (a_k-\tau(a_k)) \Big) = 0. $$
		Families of non-commutative random variable are said to be free if the $*$-subalgebras they generate are free.
		
		\item Let $ A= (a_1,\ldots ,a_k)$ be a $k$-tuple of non-commutative random variables. The joint distribution of the family $A$ is the linear form $\mu_A : P \mapsto \tau\big[ P(A, A^*) \big]$ on the set of polynomials in $2k$ non commutative indeterminates. By convergence in distribution, for a sequence of families of variables $(A_N)_{N\geq 1} = (a_{1}^{N},\ldots ,a_{k}^{N})_{N\geq 1}$ in $\mathcal C^*$-algebras $\big( \mathcal A_N, ^*, \tau_N, \norm{.} \big)$,
		we mean the pointwise convergence of the map $$ \mu_{A_N} : P \mapsto \tau_N \big[ P(A_N, A_N^*) \big],$$ and by strong convergence in distribution, we mean convergence in distribution, and pointwise convergence of the map
		$$P \mapsto \big\| P(A_N, A_N^*) \big\|.$$
		
		\item  A non commutative random variable $u$ is called a Haar unitary when it is unitary, that is $ uu^* = u^*u = 1_{\A} $, and for all $n\in \mathbb Z$, one has 
		$$
		\tau(u^n) = \left\{
		\begin{array}{ll}
		1 & \mbox{if } n=0, \\
		0 & \mbox{else.}
		\end{array}
		\right.
		$$

	\end{itemize}
	
\end{defi}

The strong convergence of non-commutative random variable is actually equivalent to the convergence of its spectrum for the Hausdorff distance. More precisely we have the following proposition whose proof can be found in \cite{collins_male} (see Proposition 2.1):

\begin{prop}
	\label{2hausdorff}
	Let $\mathbf x_N = (x_1^N,\dots, x_p^N)$ and $\mathbf x = (x_1,\dots, x_p)$ be  $p$-tuples of variables in $\mathcal C^*$-probability spaces, $(\mathcal A_N, .^*, \tau_N, \| \cdot \|)$ and $(\mathcal A, .^*, \tau, \| \cdot \|)$, with faithful states. Then, the following assertions are equivalent.
	
	\begin{itemize}
		\item $\mathbf x_N$ converges strongly in distribution to $\mathbf x$.
		\item For any self-adjoint variable $h_N=P(\mathbf x_N, \mathbf x_N^*)$, where $P$ is a fixed polynomial, $\mu_{h_N}$ converges in weak-$*$ topology to $\mu_h$ where $h=P( \mathbf x, \mathbf x^*)$. Weak-$*$ topology means relatively to continuous functions on $\mathbb C$. Moreover, the spectrum of $h_N$ converges in Hausdorff distance to the spectrum of $h$, that is, for any $\varepsilon >0$, there exists $N_0$ such that for any $N\geq N_0$, 
		\begin{equation}
		\sigma(h_N) \ \subset \ \sigma(h) \ + (-\varepsilon,\varepsilon).
		\end{equation}
	\end{itemize}
	
	\noindent In particular, the strong convergence in distribution of a single self-adjoint variable is equivalent to its convergence in distribution together with the Hausdorff convergence of its spectrum.
	
\end{prop}

It is important to note that thanks to Theorem 7.9 from \cite{nica_speicher_2006}, that we recall below, one can consider free copy of any random variable.

\begin{theorem}
	\label{2freesum}
	
	Let $(\A_i,\phi_i)_{i\in I}$ be a family of $\mathcal{C}^*$-probability spaces such that the functionals $\phi_i : \A_i\to\C$, $i\in I$, are faithful traces. Then there exist a $\mathcal{C}^*$-probability space $(\A,\phi)$ with $\phi$ a faithful trace, and a family of norm-preserving unital $*$-homomorphism $W_i: \A_i\to\A$, $i\in I$, such that:
	
	\begin{itemize}
		\item $\phi \circ W_i = \phi_i$, $\forall i \in I$.
		\item The unital $\mathcal{C}^*$-subalgebras $(W_i(\A_i))_{i\in I}$ form a free family in $(\A,\phi)$.
	\end{itemize}
\end{theorem}

\subsection{Non-commutative polynomials and derivatives}
\label{2poly}

We set $\C\langle Y_1,\dots,Y_d\rangle$ the set of non-commutative polynomials in $d$ indeterminates and in particular we fix $\PP_d = \C\langle Y_1,\dots,Y_{2d}\rangle$. Given a constant $A\in\R$, we can endow this vector space with the norm

\begin{equation}
\label{2normA}
\norm{P}_A = \sum_{M \text{ monomial}} |c_M(P)| A^{\deg M} ,
\end{equation}

\noindent where $c_M(P)$ is the coefficient of $P$ for the monomial $M$. In this subsection we define several maps on $\PP_d$ which we use multiple times in the rest of the paper, but first let us set a few notations. For $A,B,C$ non-commutative polynomials,

$$ (A\otimes B) \# C = ACB ,$$

$$ (A\otimes B) \widetilde{\#} C = BCA ,$$

$$ m(A\otimes B) = BA .$$

\noindent We define an involution $*$ on $\PP_{d}$ by fixing for all $i\in [1,d]$, $(Y_i)^* = Y_{i+d}$, $(Y_{i+d})^* = Y_i$ and then extending it to $\PP_{d}$ with the formula $(\alpha P Q)^* = \overline{\alpha} Q^* P^*$. We then define the following maps.

\begin{defi}
	\label{2application}
	
	If $1\leq i\leq d$, one set $\partial_i: \PP_{d} \longrightarrow \PP_{d} \otimes \PP_{d}$ such that for $P,Q\in \PP_{d}$,
	$$ \partial_i (PQ) = \partial_i P \times 1\otimes Q + P\otimes 1 \times \partial_i Q , $$
	$$ \partial_i Y_j = \1_{i=j} 1\otimes 1. $$
	
	\noindent We also define $D_i: \PP_{d} \longrightarrow \PP_{d}$ by $ D_i P = m\circ \partial_i P $. We similarly define $\partial_i^*$ and $D_i^*$ with the difference that for any $j$, $ \partial_i^* Y_j = \1_{i+d=j} 1\otimes 1$.
	
\end{defi}

Because they satisfy the Leibniz's rule, the maps $\partial_i$ and $\partial_i^*$ are called non-commutative derivatives. It is related to Schwinger-Dyson equations on semicircular variable, for more information see \cite{alice}, Lemma 5.4.7. While we do not use those equations in this paper, we use those associated with Haar unitary matrices. To do so, we define the following non-commutative derivative.

\begin{defi}
	\label{2application2}
	
	If $1\leq i\leq d$, one set $\delta_i: \PP_{d} \longrightarrow \PP_{d} \otimes \PP_{d}$ such that for $P,Q\in \PP_{d}$,
	$$ \delta_i (PQ) = \delta_i P \times 1\otimes Q + P\otimes 1 \times \delta_i Q , $$
	$$ \forall j\in [1,d],\quad \delta_i Y_j = \1_{i=j} Y_i\otimes 1 ,\quad  \delta_i Y_{j+d} = - \1_{i=j} 1\otimes Y_{i+d} . $$
	
	\noindent We also define $\D_i: \PP_{d} \longrightarrow \PP_{d}$ by $ \D_i P = m\circ \delta_i P $.
	
\end{defi}

We would like to apply the map $\delta_i$ to power series, more precisely the exponential of a polynomial, however since this is not well-defined in all generality we will need a few more definitions. Firstly, we need to define properly the operator norm of tensor of $\CC^*$-algebras. Since we use it later in this paper, we work with the minimal tensor product also named the spatial tensor product. For more information we refer to chapter 6 of \cite{murphy}.

\begin{defi}
	\label{2mini}
	Let $\A$ and $\B$ be $\CC^*$-algebra with faithful representations $(H_{\A},\phi_{\A})$ and $(H_{\B},\phi_{\B})$, then if $\otimes_2$ is the tensor product of Hilbert spaces, $\A\otimes_{\min}\B$ is the completion of the image of $\phi_{\A}\otimes\phi_{\B}$ in $B(H_{\A}\otimes_2 H_{\B})$ for the operator norm in this space. This definition is independent of the representations that we fixed.
\end{defi}

Consequently if $P\in \PP_d$, $ z = (z_1,\dots, z_d)$ belongs to a $\CC^*$-algebra $\A$, then $(\delta_i P^k) (z,z^*)$ belongs to $\A \otimes_{\min} \A$, and $\norm{(\delta_i P^k) (z,z^*)} \leq C_P k \norm{P(z,z^*)}^{k-1}$ for some constant $C_P$ independent of $k$. Thus we can define
\begin{equation}
\label{2extension}
(\delta_i e^{P}) (z,z^*) = \sum_{k\in \N} \frac{1}{k!} (\delta_i P^k) (z,z^*) .
\end{equation}

While we will not always be in this situation during this paper, it is important to note that if $\A = \M_N(\C)$, then up to isomorphism $\A \otimes_{\min} \A$ is simply $\M_{N^2}(\C)$ with the usual operator norm. Now we prove the following property whose proof is inspired of the one of Duhamel's formula which states that given two operators $a$ and $b$,
\begin{equation}
\label{2Duhmeqned}
e^{a}-e^{b}=\int_{0}^{1}e^{\alpha a}(a-b)e^{(1-\alpha)b}\ d\alpha.
\end{equation}

\begin{prop}
	\label{2duhamel}
	Let $P\in \PP_{d}$, $ z = (z_1,\dots, z_d)$ elements of a $\CC^*$-algebra $\A$, then
	
	$$ \left(\delta_i e^P\right)(z,z^*) = \int_0^1 \left(e^{\alpha P}\ \delta_i P\ e^{(1-\alpha)P}\right)(z,z^*) \ d\alpha , $$
	
	\noindent with convention $$A\times (B \otimes C) \times D = (AB)\otimes (CD) . $$
	
\end{prop}

\begin{proof}
	\noindent One has,
	
	\begin{align*}
	\int_0^1 \left( e^{\alpha P}\ \delta_i P\ e^{(1-\alpha)P}\right)(z,z^*) \ d\alpha &= \sum_{n,m} \int_0^1 \frac{ \alpha^n (1-\alpha)^m}{n! m!} d\alpha\ \left(P^n\ \delta_i P\ P^m\right)(z,z^*) \\
	&= \sum_k \sum_{n+m = k} \int_0^1 \frac{ \alpha^n (1-\alpha)^m}{n! m!} d\alpha\ \left(P^n\ \delta_i P\ P^m\right)(z,z^*) .\\
	\end{align*}
	
	\noindent But if $m>0$, by integration by part,
	$$ \int_0^1 \alpha^n (1-\alpha)^m d\alpha = \frac{m}{n+1}\int_0^1 \alpha^{n+1} (1-\alpha)^{m-1} d\alpha .$$
	
	\noindent Thus for any $n,m$,
	
	$$ \int_0^1 \frac{ \alpha^n (1-\alpha)^m}{n! m!} d\alpha = \int_0^1 \frac{ \alpha^{n+m}}{(n+m)!} d\alpha = \frac{1}{(m+n+1)!} .$$
	
	\noindent Hence,
	
	\begin{align*}
	\int_0^1 \left( e^{\alpha P}\ \delta_i P\ e^{(1-\alpha)P}\right)(z,z^*) \ d\alpha = \sum_k \frac{1}{(k+1)!}\sum_{n+m = k} \left(P^n\ \delta_i P\ P^m\right)(z,z^*) = \left( \delta_i\ e^P \right)(z,z^*) .
	\end{align*}
	
\end{proof}

\subsection{Free stochastic calculus}
\label{2mart}

The main idea of this paper is to use an interpolation between Haar unitary matrices and their free limit. In order to do so, we will need some notion of free stochastic calculus. The main reference in this field is the paper \cite{biane} of Biane and Speicher to which we refer for most of the proofs in this subsection. For the sake of completeness we had to introduce notations and objects that we will not necessarily use outside of this subsection. For the reader not familiar with free probability, we would suggest to focus on understanding Theorem \ref{2ito} and Definition \ref{2unitarybro}.

From now on, $(\A,\tau)$ is a $W^*$-non-commutative probability space, that is $\A$ is a von Neumann algebra, and $\tau$ is a faithful, normal (i.e. continuous for the ultraweak topology), tracial state on $\A$. We take $\A$ filtered, that is there exists a family $(\A_t)_{t\in\R^+}$ of unital, weakly closed $*$-subalgebras of $\A$,such that $\A_s\subset A_t$ for all $s\leq t$. Besides we also assume that there exist $p$ freely independent $(\A_t)_{t\in\R^+}$-free Brownian motions $((S_t^i)_{t\in\R^+})_{1\leq i\leq p}$. That is $t\mapsto S_t^i$ is weakly continuous, $S_t^i$ is a self-adjoint element of $\A_t$ with semi-circular distribution of mean $0$ and variance $t$, and for all $s\leq t$, $S_t^i-S_s^i$ is free with $\A_s$, and has semi-circular distribution of mean $0$ and variance $t-s$. Besides since the state $\tau$ is tracial, for any unital, weakly closed $*$-subalgebra $\B$ of $\A$, there exists a unique conditional expectation onto $\B$. We shall denote it by $\tau(.|\B)$. A map $t\in\R^+\mapsto M_t\in\A$ will be called a martingale with respect to the filtration $(\A_t)_{t\in\R^+}$ if for every $s\leq t$ one has $\tau(M_t|\A_s)=M_s$.

We define the opposite algebra $\A^{\op}$ as the algebra $\A$ endowed of the same addition, norm and involution, but with the product $a \times b = b\cdot a$ where $\cdot$ is the product in $\A$. We can endow $\A^{\op}$ with a faithful normal tracial state $\tau^{\op}$, which as a linear map on $\A$ is actually $\tau$. Similarly to the minimal tensor product, we will denote $L^{\infty}(\tau\otimes\tau^{\op})$ the von Neuman algebra generated by $\A\otimes\A^{\op}$ in $B(L^2(\A,\tau)\otimes_2 L^2(\A^{\op},\tau^{\op}))$ where $\otimes_2$ is the usual tensor product for Hilbert spaces. Similarly to classical stochastic calculus, we now introduce piecewise constant maps.

\begin{defi}
	A simple biprocess is a piecewise constant map $t\mapsto U_t$ from $\R^+$ into the algebraic tensor product $\A\otimes\A^{\op}$, such that $U_t = 0$ for $t$ large enough. Besides it is called adapted if for any $t\geq 0$, $U_t\in \A_t\otimes\A_t$.
\end{defi}

The space of simple biprocesses form a complex vector space that we can endow with the norm

\begin{equation}
\norm{U}_{\B^{\infty}}^2 = \int_{0}^{\infty} \norm{U_s}_{L^{\infty}(\tau\otimes\tau^{op})}^2 ds .
\end{equation}

\noindent We will denote by $\B^{\infty}_a$ the completion of the vector space of adapted simple biprocesses for this norm. Now that we have defined the notion of simple process, we can define its stochastic integral that we will later extend to  the space $\B^{\infty}_a$. 

\begin{defi}
	Let $(S_t)_{t\geq 0}$ be a free Brownian motion, $U$ be a simple adapted biprocess, we can find a decomposition $U = \sum_{j=1}^{n} A^j\otimes B^j$ and $0 = t_0\leq t_1\leq \dots\leq t_m$ such that for $t\in [t_i,t_{i+1})$, $A^j_t = A^j_{t_i}\in \A_{t_i}$ and $B^j_t = B^j_{t_i}\in \A_{t_i}^{\op}$. We define its stochastic integral by
	
	$$ \int_{0}^{\infty} U_s \# dS_s = \sum_{i=0}^{m-1} U_{t_i}\# (S_{t_{i+1}} - S_{t_i}) = \sum_{j=1}^{n} \sum_{i=0}^{m-1} A^j_{t_i} (S_{t_{i+1}} - S_{t_i}) B^j_{t_i} .$$
	
	\noindent This definition is independent of the decomposition chosen. Besides $t\mapsto \int_{0}^{t} U_s \# dX_s$ is a martingale.
	
\end{defi}

Thanks to Burkholder-Gundy inequality, that is Theorem 3.2.1 of \cite{biane}, if we see the stochastic integral as a linear map from the space of adapted simple biprocesses endowed with the norm $\norm{.}_{\B^{\infty}}$ to $\A$, then this map is continuous. Hence we can extend it to $\B^{\infty}_a$ and the martingale property remains true. Before talking about It\^o's formula, as in the classical case, we need to introduce the quadratic variation. We will not develop the idea, but by studying random matrices, in the case of simple tensors, we are prompted to define $$ \langle\langle a\otimes b , c\otimes d \rangle\rangle = a\ \tau(bc)\ d .$$

\noindent We denote by $\sharp$ the product law in $\A\otimes\A^{\op}$. If by contrast we want to use the usual product in $\A\otimes\A$, we will not put any sign. Let $\dagger$ be the linear map such that on simple tensors, $(a\otimes b)^{\dagger} = b\otimes a$. In all generality for any $Z,Y\in \A\otimes\A^{\op}$,
$$ \langle\langle Z,Y \rangle\rangle = (\1_{\A} \otimes \tau^{\op}) \left( Z \sharp (Y^{\dagger})\right) .$$

\noindent Since $ \norm{\langle\langle Z,Y \rangle\rangle} \leq \norm{Z}_{L^{\infty}(\tau\otimes\tau^{\op})} \norm{Y}_{L^{\infty}(\tau\otimes\tau^{\op})} $, we can extend this bilinear map to $Z,Y\in L^{\infty}(\tau\otimes\tau^{\op})$. Besides by Cauchy-Schwarz, for $U,V\in \B^{\infty}_a$, $\langle\langle U,V \rangle\rangle$ is integrable.

Now that we have defined all of the necessary object to do stochastic calculus, we can state It\^o's formula. We will need to handle polynomials in several processes, however Biane and Speicher only stated It\^o's formula for a product of two processes, that is if $X_0,Y_0\in\A$, $U^i,V^i\in \B_a^{\infty}$ and $K,L\in L^1(R^+,\A)$, we set

$$ Y_t = Y_0 + \int_{0}^{t} K_s ds + \sum_{1\leq i\leq p} \int_{0}^{t} U^i_s \# dS^i_s ,$$
$$ Z_t = Z_0 + \int_{0}^{t} L_s ds + \sum_{1\leq i\leq p} \int_{0}^{t} V^i_s \# dS^i_s ,$$

\noindent then for any $t\geq 0$,
\begin{align}
\label{2ito2var}
Y_t Z_t = Y_0 Z_0 &+ \int_{0}^{t} \left( Y_s L_s + K_s Z_s + \sum_{1\leq i\leq p} \langle\langle U^i_s,V^i_s\rangle\rangle\right) \ ds \\
&+ \sum_{1\leq i\leq p} \int_{0}^{t} \left((Y_s\otimes 1_{\A}) V_s^i + U_s^i (1_{\A}\otimes Z_s)\right)  \# dS^i_s . \nonumber
\end{align}

\noindent To find a proof of \eqref{2ito2var}, see Theorem 4.1.2 in \cite{biane}. While this theorem only proves the case where there is a single Brownian motion and $L=K=0$, deducing equation \eqref{2ito2var} does not require much more work. We can then deduce from equation \eqref{2ito2var}, the general It\^o's formula. Even though this formula is used without a proof by Dabrowski in \cite{dabropc}, we do not know of any satisfying reference. Hence we include a proof for self-containedness. Let us first fix a few notations.
\begin{itemize}
	\item If $P\in\C\langle X_1,\dots,X_d\rangle$, $X\in (L^{\infty}(\R^+,\A))^d$ and $K\in (L^1(\R^+,\A))^d$, then $$\partial P (X) \# K = \sum_{1\leq j\leq d} \partial_j P(X) \# K_j .$$
	\item Similarly if $U\in (\B_a^{\infty})^d$, then $\partial P(X) \sharp U = \sum_{1\leq j\leq d} \partial_j P(X) \sharp U_j$.
	\item Finally if $U,V\in \B_a^{\infty}$, $A,B,C\in L^{\infty}(\R^+,\A)$, then $(A\otimes B\otimes C)\# (U,V) = ((A\otimes B) \sharp U , (1\otimes C) \sharp V) $.
\end{itemize}

\begin{theorem}
	\label{2ito}
	
	Let $X_0\in\A^d$, $P$ be a non-commutative polynomial in $d$ indeterminates, for any $t\geq 0$, $K\in (L^1([0,t],\A))^d$ and for every $i\in[1,p]$, $(\1_{s\leq t} U^i_s)_{s\in\R^+} \in (\B_a^{\infty})^d$. With $I$ the identity map on $\PP_d$, we define 
	$$ X_t = X_0 + \int_{0}^{t} K_s ds + \sum_{1\leq i\leq p} \int_{0}^{t} U^i_s \# dS^i_s ,$$
	$$ \Delta_U(P) (X)= \sum_{1\leq i\leq p} \sum_{1\leq j,k\leq d} \langle\langle\  \left((\partial_j \otimes I)\circ \partial_k P(X)\right) \# (U^{i,j},U^{i,k}) \  \rangle\rangle .$$
	\noindent Then for any $t\geq 0$, $\partial P (X) \# K$ and $\Delta_U(P) (X) \in L^1([0,t],\A))$, and $(\1_{s\leq t} \partial P(X_s) \sharp U_s) _{s\in\R^+} \in \B_a^{\infty}$. Finally for any $t\geq 0$,
	$$ P(X_t) = P(X_0) + \int_0^t \partial P (X_s) \# K_s\  ds + \sum_{1\leq i\leq p} \int_0^t \left(\partial P(X) \sharp U^i_s \right)  \#dS_s^i + \int_{0}^t \Delta_U(P) (X_s)\  ds .$$
\end{theorem}

\begin{proof}
	Thanks to Burkholder-Gundy inequality, that is Theorem 3.2.1 of \cite{biane}, we know that
	$$ \sup_{0\leq s\leq t} \norm{X_s^j} \leq \norm{X_0^j} + \norm{K_s^j}_{L^1([0,t],\A)} + \sum_{1\leq i\leq p} \norm{U^{i,j} \1_{[0,t]}}_{\B^{\infty}_a} .$$
	
	\noindent Thus for any $t\in\R^+$, $(X_s)_{s\in [0,t]} \in L^{\infty}([0,t],\A)^d$, hence for any polynomial $P$, $\partial P (X) \# K \in L^{\infty}([0,t],\A)$, and $(\1_{s\leq t} \partial P(X_s) \sharp U_s) _{s\in\R^+} \in \B_a^{\infty}$. Given that 
	$$ \norm{\langle\langle Z,Y \rangle\rangle} \leq \norm{Z}_{L^{\infty}(\tau\otimes\tau^{\op})} \norm{Y}_{L^{\infty}(\tau\otimes\tau^{\op})}, $$
	we also have that $\Delta_U(P) (X) \in L^1([0,t],\A))$. Finally to prove the formula, we proceed recurrently. If $P$ is of degree $1$, there is nothing to prove. For larger degree, by linearity we only need to deal with the case where $P$ is a monomial. Thus we can write $P = QR$ with $Q$ and $R$ monomials of smaller degree for which the formula is verified. Thus thanks to equation \eqref{2ito2var}, we have that
	
	\begin{align*}
	P(X_t) =&\ Q(X_0)R(X_0) + \int_{0}^{t} Q(X_s)\times (\partial R(X_s) \# K_s) + (\partial Q(X_s)\# K_s)\times  R(X_s) ds \\
	&+ \int_{0}^{t} Q(X_s)\ \Delta_U(R) (X_s)\  + \Delta_U(Q) (X_s)\ R(X_s)\  + \sum_{1\leq i\leq p} \langle\langle\ \partial Q(X_s) \sharp U^i_s, \partial R(X_s) \sharp U^i_s\ \rangle\rangle\  ds \\
	&+ \sum_{1\leq i\leq p} \int_{0}^{t} \left(\ (Q(X_s)\otimes \1_{\A})\times (\partial R(X_s) \sharp U^i_s) + (\partial Q(X_s) \sharp U^i_s)\times  (\1_{\A}\otimes R(X_s))\ \right) \#dS_s^i .
	\end{align*}
	
	\noindent It is clear that,
	$$ \partial (QR)(X_s) \# K_s = Q(X_s)\times (\partial R(X_s) \# K_s) + (\partial Q(X_s)\# K_s)\times  R(X_s) ,$$
	$$ \partial (QR)(X_s) \sharp U^i_s = (Q(X_s)\otimes \1_{\A})\times (\partial R(X_s) \sharp U^i_s) + (\partial Q(X_s) \sharp U^i_s)\times  (\1_{\A}\otimes R(X_s)) . $$
	
	\noindent And finally,
	\begin{align*}
	\Delta_U(QR) (X) =& \sum_{1\leq i\leq p} \sum_{1\leq j,k\leq d} \langle\langle\  \left(\ (\partial_j \otimes I)\circ \partial_k (QR)(X)\ \right) \# (U^{i,j},U^{i,k}) \  \rangle\rangle \\
	=&  \sum_{1\leq i\leq p} \sum_{1\leq j,k\leq d} \sum_{Q = A X_j B X_k C}  \langle\langle\ (A(X)\otimes B(X)\otimes C(X)) \# (U^{i,j},U^{i,k}) \  \rangle\rangle R(X) \\
	&\quad\quad\quad\quad\quad + \sum_{R = A X_j B X_k C}  Q(X) \langle\langle\ (A(X)\otimes B(X)\otimes C(X)) \# (U^{i,j},U^{i,k}) \  \rangle\rangle \\
	&\quad\quad\quad\quad\quad + \sum_{Q = A X_j B, R= C X_k D}  \langle\langle\ (A(X)\otimes (BC)(X)\otimes D(X)) \# (U^{i,j},U^{i,k}) \  \rangle\rangle \\
	=& Q(X)\ \Delta_U(R) (X)\  + \Delta_U(Q) (X)\ R(X)\  + \sum_{1\leq i\leq p} \langle\langle\ \partial Q(X) \sharp U^i, \partial R(X) \sharp U^i\ \rangle\rangle .
	\end{align*}
	
\end{proof}

Finally, one of the fundamental tool that we use in this paper is the free unitary Brownian motion, a good reference on the matter is \cite{freebrlaw}. In particular one can find a proof of its existence.

\begin{defi}
	\label{2unitarybro}
	Let $(S_t)_{t\geq 0}$ be a free Brownian motion adapted to a filtered $W^*$-probability space $(\A,(\A_t)_{t\geq 0},\tau)$, the free unitary Brownian motion $(u_t)_{t\geq 0}$ is the unique solution to the equation
	\begin{equation}
	\label{2defunbro}
	\forall t\geq 0, \quad  u_t = 1_{\A} - \int_0^t \frac{u_s}{2}\ ds + \i \int_0^t (u_s\otimes 1_{\A}) \# dS_s .
	\end{equation}
	
	\noindent In particular, for any $t\geq 0$, $u_t$ is unitary, that is $u_t u_t^* = u_t^* u_t = 1_{\A}$.
\end{defi}

Although we do not use this notation in this paper, similarly to the classical case, it is usual to write equation \eqref{2defunbro} as
$$ u_0 = 1_{\A}, \quad  du_t =  - \frac{u_t}{2}\ dt + \i (u_t\otimes 1_{\A}) \# dS_t . $$

\subsection{Notations}
\label{2notationbasic}

Let us now fix a few notations concerning the spaces and traces that we use in this paper.

\begin{defi}
	\label{2bases}
	\begin{itemize}
		\item $(\A_N,\tau_N)$ is the free product of $\M_N(\C)$ with the von Neuman algebra $\A$ from the former subsection. To build $\A_N$ we use Theorem \ref{2freesum} and we get a $\mathcal{C}^*$-probability space $C$ with a faithful trace $\varphi$. Since we want $(\A_N,\tau_N)$ to be a von Neuman algebra, we set $L^2(C,\varphi)$ as the completion of $C$ for the norm $a\mapsto \phi(a^*a)^{1/2}$, we have an injective $\CC^*$-algebra morphism from $C$ to $B(L^2(C,\varphi))$. We then proceed to take $\A_N$ the closure of the image of $C$ in this space for the weak topology. As for $\tau_N$, since we can extend $(x,y) \mapsto \varphi(x^* y)$ to a scalar product $\langle .,. \rangle_{\varphi}$ on $L^2(C,\varphi)$, we set for $a\in B(L^2(C,\varphi))$, $\tau_N(a) = \langle a(1),1 \rangle_{\varphi}$.
		\item Note that when restricted to $\M_N(\C)$, $\tau_N$ is just the regular renormalized trace on matrices. Similarly we will usually denote $\tau_M$ and $\tau_k$ the renormalized trace on $M_M(\C)$ and $M_k(\C)$. As in the former subsection, the restriction of $\tau_{N}$ to the $\mathcal{C}^*$-algebra $\A$ is denoted as $\tau$. 
		\item $\tr_N$ is the non-renormalized trace on $\M_N(\C)$.
		\item $E_{i,j}$ is the matrix whose only non-zero coefficient is $(i,j)$ and this coefficient has value $1$, the size of the matrix $E_{i,j}$ will depend on the context.
		\item In general we identify $\M_N(\C)\otimes \M_k(\C)$ with $\M_{kN}(\C)$ through the isomorphism $ E_{i,j}\otimes E_{r,s} \mapsto E_{i+(r-1)N,j+(s-1)N} $, similarly we identify $\tr_N\otimes\tr_k$ with $\tr_{kN}$.
		\item If $A^N=(A_1^N,\dots,A_d^N)$ and $B^M=(B_1^M,\dots,B_d^M)$ are two families of matrices, then we denote $A^N\otimes B^M=(A_1^N\otimes B^M_1,\dots,A_d^N\otimes B^M_d)$ and if $M=N$, $A^N B^N=(A_1^N B^N_1,\dots,A_d^N B^N_d)$. We typically use the notation $X^N\otimes I_M$ for the family $(X^N_1\otimes I_M,\dots,X^N_d\otimes I_M)$.
		\item If $P\in\PP_d$, in order to avoid cumbersome notations when evaluating $P$ in $(X,X^*)$, instead of denoting $P(X,X^*)$ we will write $\widetilde{P}(X)$.
		\item We define $(e_i)_{i\in [1,M]}$, $(g_i)_{i\in [1,N]}$ and $(f_i)_{i\in [1,k]}$ the canonical basis of $\C^M$, $\C^N$ and $\C^k$
	\end{itemize}
\end{defi}

A polynomial $P\in \PP_{d}$ is said to be self-adjoint if $P^* = P$. Self-adjoint polynomials have the property that if $z_1,\dots,z_d$ are elements of a $\mathcal{C}^*$-algebra, then  $P(z_1,\dots,z_d,z_1^*,\dots,z_d^*)$ is self-adjoint. Now that we have defined the notion of self-adjoint polynomial we give a property which justifies computations that we will do later on:

\begin{prop}
	\label{2mesurable}
	Let the following objects be given,
	\begin{itemize}
		\item $u=(u_t^1,\dots ,u_t^p)_{t\geq 0}$ a family of $p$ freely independent free unitary Brownian motions,
		\item $f\in \mathcal{C}(\R) $ the set of continuous function on $\R$,
		\item $P$ a self-adjoint polynomial.
	\end{itemize}
	Then with $\U_N$ the group of unitary matrices of size $N$, the following map is measurable:
	$$ (U^N,Z^{NM})\in \U_N^p \times \M_{NM}(\C)^{d-p}  \mapsto \tau_N\otimes\tau_M\left(f\left(\widetilde{P}\left((U^N u_t)\otimes I_M,Z^{NM}\right)\right)\right) .$$
\end{prop}

For a full proof we refer to \cite{un}, Proposition 2.7. But in a few words, it is easy to see the measurability when $f$ is a polynomial since then this map is also polynomial in the coefficient of $U^N$ and $Z^{NM}$, and we conclude by density. Actually we could easily prove that this map is continuous, however we do not need it. The only reason we need this property is to justify that if $U^N$ is a family of random matrices, then the random variable $\tau_N\otimes\tau_M\left(f\left(\widetilde{P}(u_t^N\otimes I_M,Z^{NM})\right)\right)$ is well-defined. To conclude this subsection we introduce different notations related to maps defined on tensor spaces.

\begin{defi}
	Let $n : A\otimes B \in\M_M(\C)^{\otimes 2}\mapsto AB \in \M_M(\C)$, we define the linear map $(\tau_N\otimes I_M) \bigotimes (\tau_N\otimes I_M) : (\A_N\otimes \M_M(\C))^{\otimes 2} \to M_M(\C)$ as
	$$(\tau_N\otimes I_M) \bigotimes (\tau_N\otimes I_M) = n\circ (\tau_N\otimes I_M)^{\otimes 2} .$$
	We will also use the shorter notation $(\tau_N\otimes I_M)^{\bigotimes 2}$.
\end{defi}

\subsection{Random matrix model}

We conclude this section by giving the definition and a few properties on the models of random matrices that we will study. 

\begin{defi}
	\label{2Haardef}
	A Haar unitary matrix of size $N$ is a random matrix distributed according to the Haar measure on the group of unitary matrices of size $N$.
\end{defi}

\begin{defi}
	\label{2HBdef}
	A Hermitian Brownian motion $(X_t^N)_{t\in\R^+}$ of size $N$ is a self-adjoint matrix whose coefficients are random variables with the following laws:
	\begin{itemize}
		\item For $1\leq i\leq N$, the random variables $\sqrt{N} ((X^N_t)_{i,i})_{t\in\R^+}$ are independent Brownian motions.
		\item For $1\leq i<j\leq N$, the random variables $(\sqrt{2N}\ \Re{(X^N_t)_{i,j}})_{t\in\R^+}$ and $(\sqrt{2N}\ \Im{(X^N_t)_{i,j}})_{t\in\R^+}$ are independent Brownian motions, independent of $\sqrt{N} ((X^N_t)_{i,i})_{t\in\R^+}$.
	\end{itemize}
\end{defi}

To study the free unitary Brownian motion, we will need to study its finite dimensional equivalent, the unitary Brownian motion. Typically it is defined as the Markov process whose infinitesimal generator is the Laplacian operator on the unitary group. However given the upcoming computations in this paper, it is better to use an equivalent definition as the solution of a stochastic differential equation. We refer to subsection 2.1 of \cite{UBMN} for a short summary on the different definitions.

\begin{defi}
	\label{2UBdef}
	Let $X^N$ be a Hermitian Brownian motion, then the unitary Brownian motion $(U^N_t)_{t\geq 0}$ is the solution of the following stochastic differential equation:
	\begin{equation}
	\label{2defbrunit}
	dU_t^N = \i U_t^N dX_t^N - \frac{1}{2} U_t^N dt,\quad U_0^N = I_N,
	\end{equation}
	
	\noindent where we formally define $U_t^N dX_t^N$ by simply taking the matrix product $$(U_t^N dX_t^N)_{i,j} = \sum_k (U_t^N)_{i,k} d(X_t^N)_{k,j}.$$
	
	\noindent In particular, almost surely, for any $t$, $U_t^N$ is a unitary matrix of size $N$.
\end{defi}

The following property is typical of the kind of computation that we can do with unitary Brownian motion with classical stochastic calculus, see \cite{UBMN} for example.

\begin{prop}
	\label{2unitarybr}
	
	Let $U^N_1,\dots,U^N_p$ be independent unitary Brownian motions of size $N$, $A^N_{p+1},\dots,A^N_d$ be deterministic matrices, $Q\in\PP_d$ be a monomial, we set $Q_s$ the monomial evaluated in $(U^N_{1,s},\dots,U^N_{p,s},$ $A^N_{p+1},\dots,A^N_d)$ and their adjoints, $|Q|_B$ the degree of $Q$ with respect to $(U_1,\dots,U_p,U_1^*,\dots,U_p^*)$. Then there exists a martingale $J$ such that,
	
	\begin{align*}
	d\tr_N Q_s = dJ_s - \frac{|Q|_B}{2} \tr_N Q_s\ ds &- \frac{1}{N} \sum_{i\leq p,\  Q = A U_i B U_i C} \tr_N(A_s U^N_{i,s} C_s) \tr_N(B_s U^N_{i,s})\ ds\\
	&- \frac{1}{N} \sum_{i\leq p,\  Q = A U_i^* B U_i^* C} \tr_N\left(A_s {U^N_{i,s}}^* C_s\right) \tr_N\left(B_s {U^N_{i,s}}^*\right)\ ds \\
	&+ \frac{1}{N} \sum_{i\leq p,\ Q=A U_i B U_i^* C } \tr(A_s C_s) \tr(B_s)\ ds \\
	&+ \frac{1}{N} \sum_{i\leq p,\ Q=A U_i^* B U_i C } \tr(A_s C_s) \tr(B_s)\ ds .
	\end{align*}
	
\end{prop}

%
%
%

\section{Preliminaries}
\label{2preli}

\subsection{A matricial inequality}

\noindent We are indebted to Mikael de la Salle for supplying us with the proof of the following lemma.

\begin{lemma}[de la Salle]
	\label{2normineq}
	Let $\A$ be a $\mathcal{C}^*$-algebra, $A_1,A_2\in \A$, $B_1,B_2\in \M_M(\C)$, as in subsection \ref{2mart} we define $(A_1\otimes B_1)\sharp (A_2\otimes B_2) = (A_1A_2)\otimes(B_2B_1)$. Then if $x,y\in \A\otimes \M_M(\C)$, with the operator norm in their respective space,
	$$ \norm{x\sharp y} \leq M \norm{x}\norm{y} .$$
\end{lemma}

\begin{proof}
	
	We write $x=\sum_{ 1\leq i,j \leq M} x_{i,j}\otimes E_{i,j}$, $y = \sum_{ 1\leq k,l \leq M} y_{k,l}\otimes E_{k,l}$, then 
	$$x\sharp y = \sum_{i,j,k} x_{k,j} y_{i,k} \otimes E_{i,j} . $$
	
	\noindent We define $A_k = \sum_{i,j} x_{k,j} y_{i,k} \otimes E_{i,j}$, $X_k = \sum_j x_{k,j}\otimes E_{k,j}\otimes I_M$, $Y_k = \sum_i y_{i,k}\otimes I_M \otimes E_{i,k}$. Then by using the fact that $X_k$ and $Y_k$ are band matrices, we have $\norm{X_k}\leq \norm{x}$ and $\norm{Y_k}\leq \norm{y}$. Besides $\norm{x\sharp y} \leq \sum_{1\leq k\leq M} \norm{A_k}$. Finally	we have for any $k$,
	
	\begin{align*}
	\norm{X_kY_k}^2 &= \norm{X_kY_k (X_kY_k)^*} \\
	&= \norm{\left(\sum_{i,j} x_{k,j} y_{i,k} \otimes E_{k,j}\otimes E_{i,k}\right)\left(\sum_{i,j} x_{k,j} y_{i,k} \otimes E_{k,j}\otimes E_{i,k}\right)^*} \\
	&= \norm{ \sum_{i,j,u,v} x_{k,j} y_{i,k} y_{u,k}^* x_{k,v}^* \otimes E_{k,j}E_{v,k} \otimes E_{i,k}E_{k,u} } \\
	&= \norm{ \sum_{i,j,u} x_{k,j} y_{i,k} y_{u,k}^* x_{k,v}^* \otimes E_{k,k} \otimes E_{i,j}E_{v,u}} \\
	&= \norm{A_kA_k^*\otimes E_{k,k}} \\
	&= \norm{A_k}^2 .
	\end{align*}
	
	\noindent Thus  $\norm{x\sharp y} \leq \sum_{1\leq k\leq M} \norm{A_k}\leq M \norm{x}\norm{y}$.
	
\end{proof}

\subsection{A Poincar\'e type equality}

One of the main tool when dealing with GUE random matrices is the Poincar\'e inequality (see Definition 4.4.2 from \cite{alice}), which gives us a sharp majoration of the variance of a function in these matrices. Typically this inequality shows that the variance of the renormalized trace of a polynomial in GUE random matrices, which a priori is of order $\mathcal{O}(1)$, is of order $\mathcal{O}(N^{-2})$. In this paper we need a similar type of inequality but instead of working with independent GUE random matrices, we work with marginal of independent unitary Brownian motions at times $t$. We will follow the approach of \cite{levy-maida}, Proposition 6.1. I would like to thank one of the reviewer for pointing out that an alternative approach could be to use the more general Theorem 4.3 from \cite{rapport}. This theorem deals with Brownian motions defined with the help of two parameters $r$ and $s$, and the unitary Brownian motion matches with the case where $(r,s) = (1,0)$. However the proof of Proposition \ref{2concentration} that we give below is a good introduction to the kind of computations that we deal with in this paper unlike the approach taken in \cite{rapport} and needs less prerequisite to be understood.

\begin{prop}
	\label{2concentration}
	
	Let $Q\in\PP_d$, $(U_t^N)_{t\in\R^+}$, $(V_t^N)_{t\in\R^+}$, $(W_t^N)_{t\in\R^+}$ be independent families of $p$ unitary Brownian motions of size $N$. Let $A^N$ be a family of $d-p$ deterministic matrices, with notations as in Definition \ref{2bases}, one has for any $T\geq 0$,
	
	\begin{equation*}
	\var\Big(\tr_N\left(\widetilde{Q}(U_T^N,A^N)\right)\Big) = \frac{1}{N} \sum_{1\leq k\leq p} \int_{0}^{T} \E\Big[ \tr_N\Big( \widetilde{\D_kQ}(V_{T-t}^N U_t^N, A^N) \times \widetilde{\D_kQ}(W_{T-t}^N U_t^N, A^N)^* \Big) \Big] dt .
	\end{equation*}
	
\end{prop}

\begin{proof}
	To simplify notations, we will not write the index $N$ in $U_t^N,V_t^N,W_t^N$ and $A^N$. For $U\in \M_N(\C)^p$, we set $f:(U,U^*)\mapsto \tr_N(Q(U,A,U^*,A^*))$. We can view $f$ as a polynomial in the coefficients of the matrices $U$ and their conjuguates, since those are complex variables we use the notion of complex differential. That is if $g: (x,y)\in\R^2\to g(x,y)\in\C$ is a differentiable function, we define $\partial_z g = \frac{1}{2} \left(\partial_x g - \i \partial_y g\right)$ and $\partial_{\overline{z}} g = \frac{1}{2} \left(\partial_x g + \i \partial_y g\right)$. If $u_k^{i,j}$ is the $(i,j)$-coefficient of the $k$-th matrix in $U$, we denote the differential of $f$ with respect to $u_k^{i,j}$ by $\partial_{u_k^{i,j}}f$, and the differential of $f$ with respect to the conjuguate of this coefficient by $\partial^*_{u_k^{j,i}}f$. In particular,
	\begin{align*}
	&\partial_{u_k^{i,j}}( (U_k)_{i,j}) = 1,\quad \partial_{u_k^{i,j}}^*( (U_k^*)_{i,j}) = 1, \\
	&\partial_{u_k^{i,j}}( (U_k)_{a,b}) = 0,\quad \partial_{u_k^{i,j}}^*( (U_k^*)_{a,b}) = 0, \text{ for all } (a,b)\neq (i,j), \\
	&\partial_{u_k^{i,j}}( (U_k^*)_{a,b}) = 0,\quad \partial_{u_k^{i,j}}^*( (U_k)_{a,b}) = 0, \text{ for any } (a,b).
	\end{align*}
	
	\noindent Next we introduce
	$$ M_t = P_{T-t}f (U_t,U_t^*),$$
	
	\noindent where $P_{T-t} f(U,U^*) = \E_V[f(V_{T-t}U,(V_{T-t}U)^*)]$ with $ (V_t)_{t\geq 0} $, $p$ independent unitary Brownian motions of size $N$ and $E_V$ the expectation with respect to $ (V_t)_{t\geq 0} $. We will follow the approach of \cite{levy-maida}, Proposition 6.1, and show that $(M_t)_{0\leq t\leq T}$ is a martingale. It will follow that,
	\begin{align}
	\label{2explanation}
	\var\Big(\tr_N\left(\widetilde{Q}(U_T^N,A^N)\right)\Big) &= \E[ |f(U_T,U_T^*)|^2 - |\E[f(U_T,U_T^*)]|^2 ] \nonumber \\
	&= \E[ M_T \overline{M_T} - M_0\overline{M_0} ] \\
	&= \E\left[ \langle M_T, \overline{M_T}\rangle \right] . \nonumber
	\end{align}
	
	\noindent If we set $ (X_t)_{t\geq 0} $, $p$ independent Hermitian Brownian motions of size $N$, and $f_t = P_{T-t} f$, then
	
	\begin{align*}
	dM_t = (\partial_tf_t) (U_t,U_t^*) dt &+ \sum_{1\leq k\leq p,\ 1\leq i,j\leq N } (\partial_{u_{k}^{i,j}} f_t)(U_t,U_t^*)\ d(U_{k,t})_{i,j} + (\partial_{u_{k}^{i,j}}^* f_t)(U_t,U_t^*)\ d(U_{k,t}^*)_{i,j} \\
	&+ \frac{1}{2} \sum_{\substack{1\leq k\leq p,\ 1\leq i,j,r,s\leq N \\ \varepsilon_1,\varepsilon_2\in \{1,*\} }} (\partial_{u_{k}^{i,j}}^{\varepsilon_1} \partial_{u_{k}^{r,s}}^{\varepsilon_2} f_t)(U_t,U_t^*)\ d\langle (U_{k,t}^{\varepsilon_1})_{i,j} , (U_{k,t}^{\varepsilon_2})_{r,s} \rangle_t .
	\end{align*}
	
	\noindent By using equation \eqref{2defbrunit}, we can isolate the martingale term in the previous equation. As for the term associated to $dt$, as long as we show that $(M_t)_{t\geq 0}$ is a martingale it will be $0$. To do so, we set $\mathcal F_t = \sigma((U_s)_{0\leq s\leq t})$, then $\mathcal F_t$ is a filtration adapted to $(M_t)_{t\geq 0}$. Besides if we set $N_t = \E[f (U_T,U_T^*)\ |\ \mathcal F_t]$, then since if we redefine $V_{T-t} = U_TU_t^*$, it is still a family of $p$ independent unitary Brownian motions of size $N$, independent of $\mathcal F_t$. This implies that 
	$$ N_t = \E[f (U_T,U_T^*)\ |\ \mathcal F_t] = \E[f (V_{T-t}U_t,(V_{T-t}U_t)^*)\ |\ \mathcal F_t] = \E_V[f (V_{T-t}U_t,(V_{T-t}U_t)^*)] = M_t .$$
	
	\noindent Hence $(M_t)_{t\geq 0}$ is a martingale and 	
	\begin{equation}
	\label{2martingalebra}
	dM_t = \i \sum_{1\leq k\leq p,\ 1\leq i,j\leq N} (\partial_{u_{k}^{i,j}} f_t)(U_t,U_t^*)\ (U_{k,t}\ dX_{k,t}^N)_{i,j} - (\partial_{u_{k}^{i,j}}^* f_t)(U_t,U_t^*)\ ((U_{k,t}\ dX_{k,t}^N)^*)_{i,j} ,
	\end{equation}
	
	\noindent Therefore, as we saw in equation \eqref{2explanation}, we are left with computing the bracket of $M_t$. To begin with we have,
	$$ \langle (U_{k,t}\ dX_{k,t}^N)_{i,j}, \overline{(U_{k,t}\ dX_{k,t}^N)_{r,s}} \rangle = \1_{i=r,j=s} \frac{dt}{N} , $$
	$$ \langle ((U_{k,t}\ dX_{k,t}^N)^*)_{i,j}, \overline{((U_{k,t}\ dX_{k,t}^N)^*)_{r,s}} \rangle = \1_{i=r,j=s} \frac{dt}{N} , $$
	$$ \langle (U_{k,t}\ dX_{k,t}^N)_{i,j}, \overline{((U_{k,t}\ dX_{k,t}^N)^*)_{r,s}} \rangle = (U_{k,t})_{i,r} (U_{k,t})_{s,j} \frac{dt}{N}  , $$
	$$ \langle ((U_{k,t}\ dX_{k,t}^N)^*)_{i,j}, \overline{(U_{k,t}\ dX_{k,t}^N)_{r,s}} \rangle = (U_{k,t}^*)_{s,j} (U_{k,t}^*)_{i,r} \frac{dt}{N} .$$	
	\break
	
	\noindent Let us remind that $f:(U,U^*)\mapsto \tr_N(Q(U,A,U^*,A^*))$, thus one has
	
	$$ (\partial_{u_{k}^{i,j}} f_t)(U_t,U_t^*) = \E_V\left[ \tr_N(\widetilde{D_kQ}(V_{T-t}U_t,A)\ V_{k,T-t} E_{i,j}) \right] ,$$
	
	$$ (\partial_{u_{k}^{i,j}}^* f_t)(U_t,U_t^*) = \E_V\left[ \tr_N(V_{k,T-t}^*\ \widetilde{D_k^*Q}(V_{T-t}U_t,A)\ E_{i,j}) \right] .$$
	
	\noindent We will now compute four different brackets, and by summing them we will get the bracket of $M_t$ (see equation \eqref{2martingalebra}). First,
	\begin{align}
	\label{2cas11}
	&\left< \sum_{i,j,k} (\partial_{u_{k}^{i,j}} f_t)(U_t,U_t^*)\ (U_{k,t}\ dX_{k,t}^{N})_{i,j}, \overline{\sum_{i,j,k} (\partial_{u_{k}^{i,j}} f_t)(U_t,U_t^*)\ (U_{k,t}\ dX_{k,t}^N)_{i,j}} \right> \\
	&= \sum_k \sum_{i,j,r,s} (\partial_{u_{k}^{i,j}} f_t)(U_t,U_t^*)\ \overline{(\partial_{u_{k}^{r,s}} f_t)(U_t,U_t^*)} \left< (U_{k,t}\ dX_{k,t}^{N})_{i,j}, \overline{(U_{k,t}\ dX_{k,t}^{N})_{r,s}} \right> \nonumber \\
	&= \frac{1}{N} \sum_k \sum_{i,j} (\partial_{u_{k}^{i,j}} f_t)(U_t,U_t^*)\ \overline{(\partial_{u_{k}^{i,j}} f_t)(U_t,U_t^*)} \ dt \nonumber \\
	&= \frac{1}{N} \sum_k \sum_{i,j} \E_V\left[ \tr_N(\widetilde{D_kQ}(V_{T-t}U_t,A)\ V_{k,T-t} E_{i,j}) \right] \E_V\left[ \tr_N(E_{j,i} (V_{k,T-t})^* \widetilde{D_kQ}(V_{T-t}U_t,A)^*\  ) \right] \ dt \nonumber \\
	&= \frac{1}{N} \sum_k \E_{V,W}\left[ \tr_N(\widetilde{D_kQ}(V_{T-t}U_t,A)\ V_{k,T-t} W_{k,T-t}^* \widetilde{D_kQ}(W_{T-t}U_t,A)^* ) \right] dt . \nonumber 
	\end{align}
	
	\noindent Similarly one has,
	
	\begin{align}
	\label{2cas12}
	&\left< \sum_{i,j,k} (\partial_{u_{k}^{i,j}}^* f_t)(U_t,U_t^*)\ ((U_{k,t}\ dX_{k,t}^{N})^*)_{i,j}, \overline{\sum_{i,j,k} (\partial_{u_{k}^{i,j}}^* f_t)(U_t,U_t^*)\ ((U_{k,t}\ dX_{k,t}^{N})^*)_{i,j}} \right> \\
	&= \frac{1}{N} \sum_k \E_{V,W}\left[ \tr_N\left( V_{k,T-t}^*\ \widetilde{D_k^*Q}(V_{T-t}U_t,A)\ \widetilde{D_k^*Q}(W_{T-t}U_t,A)^*\ W_{k,T-t} \right) \right] dt . \nonumber 
	\end{align}
	
	\noindent Next we have,
	
	\begin{align}
	\label{2cas13}
	&\left< \sum_{i,j,k} (\partial_{u_{k}^{i,j}} f_t)(U_t,U_t^*)\ (U_{k,t}\ dX_{k,t}^{N})_{i,j}, \overline{\sum_{i,j,k} (\partial_{u_{k}^{i,j}}^* f_t)(U_t,U_t^*)\ ((U_{k,t}\ dX_{k,t}^{N})^*)_{i,j}} \right> \\
	&= \sum_k \sum_{i,j,r,s} (\partial_{u_{k}^{i,j}} f_t)(U_t,U_t^*)\ \overline{(\partial_{u_{k}^{r,s}}^* f_t)(U_t,U_t^*)} \left< (U_{k,t}\ dX_{k,t}^{N})_{i,j}, \overline{((U_{k,t}\ dX_{k,t}^{N})^*)_{r,s}} \right> \nonumber \\
	&= \frac{1}{N} \sum_k \sum_{i,j,r,s} (\partial_{u_{k}^{i,j}} f_t)(U_t,U_t^*)\ \overline{(\partial_{u_{k}^{r,s}}^* f_t)(U_t,U_t^*)}\ (U_{k,t})_{i,r} (U_{k,t})_{s,j}\ dt \nonumber \\
	&= \frac{1}{N} \sum_k \sum_{i,j,r,s} \E_V\left[ \tr_N(\widetilde{D_kQ}(V_{T-t}U_t,A)\ V_{k,T-t} E_{i,j}) \right] \nonumber \\
	&\quad\quad\quad\quad\quad\quad \times \E_W\left[ \tr_N(E_{s,r} \widetilde{D_k^*Q}(W_{T-t}U_t,A)^*\ W_{k,T-t} ) \right] \ (U_{k,t})_{i,r} (U_{k,t})_{s,j}\ dt \nonumber \\
	&= \frac{1}{N} \sum_k \E_{V,W}\Bigg[ \sum_{i,j,r,s} \left( \widetilde{D_kQ}(V_{T-t}U_t,A)\ V_{k,T-t} \right)_{j,i} (U_{k,t})_{i,r} \\
	&\quad\quad\quad\quad\quad\quad\quad\quad\quad\times \left(\widetilde{D_k^*Q}(W_{T-t}U_t,A)^*\ W_{k,T-t}\right)_{r,s} (U_{k,t})_{s,j}  \Bigg] dt \nonumber \\
	&= \frac{1}{N} \sum_k \E_{V,W}\left[ \tr_N\left( \widetilde{D_kQ}(V_{T-t}U_t,A)\ V_{k,T-t} U_{k,t}\ \widetilde{D_k^*Q}(W_{T-t}U_t,A)^*\ W_{k,T-t} U_{k,t} \right) \right] dt . \nonumber 
	\end{align}
	
	\noindent And similarly,
	
	\begin{align}
	\label{2cas14}
	&\left< \sum_{i,j,k} (\partial_{u_{k}^{i,j}}^* f_t)(U_t,U_t^*)\ ((U_{k,t}\ dX_{k,t}^{N})_{i,j}, \overline{\sum_{i,j,k} (\partial_{u_{k}^{i,j}} f_t)(U_t,U_t^*)\ (U_{k,t}\ dX_{k,t}^{N})_{i,j}} \right> \\
	&= \frac{1}{N} \sum_k \E_{V,W}\left[ \tr_N\left( (V_{k,T-t} U^k_t)^* \widetilde{D_k^*Q}(V_{T-t}U_t,A)\ (W_{k,T-t} U^k_t)^* \widetilde{D_kQ}(W_{T-t}U_t,A)^* \right) \right] dt . \nonumber 
	\end{align}
	
	\noindent We sum equations \eqref{2cas11} to \eqref{2cas14}.
	\begin{align*}
	\var\Big(\tr_N(\widetilde{Q}(U_T^N,A^N))\Big)& \\
	= \frac{1}{N} \sum_k \int_{0}^{T} \E\Big[ \tr_N\Big( &\widetilde{D_kQ}(V_{T-t}U_t,A)\ V_{k,T-t}U_{k,t}\ (W_{k,T-t}U_{k,t})^* \widetilde{D_kQ}(W_{T-t}U_t,A)^*  \\
	&+ (V_{k,T-t}U_{k,t})^*\ \widetilde{D_k^*Q}(V_{T-t}U_t,A)\ \widetilde{D_k^*Q}(W_{T-t}U_t,A)^*\ W_{k,T-t}U_{k,t} \\
	&- \widetilde{D_kQ}(V_{T-t}U_t,A)\ V_{k,T-t} U_{k,t}\ \widetilde{D_k^*Q}(W_{T-t}U_t,A)^*\ W_{k,T-t} U_{k,t} \\
	&- (V_{k,T-t} U^k_t)^* \widetilde{D_k^*Q}(V_{T-t}U_t,A)\ (W_{k,T-t} U^k_t)^* \widetilde{D_kQ}(W_{T-t}U_t,A)^* \Big) \Big] dt \\
	= \frac{1}{N} \sum_k \int_{0}^{T} \E\Big[ \tr_N\Big( &\left(\widetilde{D_kQ}(V_{T-t}U_t,A)\ V_{k,T-t}U_{k,t} - (V_{k,T-t}U_{k,t})^*\ \widetilde{D_k^*Q}(V_{T-t}U_t,A) \right) \times \\
	& \left( \widetilde{D_kQ}(W_{T-t}U_t,A) W_{k,T-t} U_{k,t} -  (W_{k,T-t} U_{k,t})^* \widetilde{D_k^*Q}(W_{T-t}U_t,A) \right)^* \Big) \Big] dt .
	\end{align*}
	
	\noindent Hence the conclusion.
	
\end{proof}

\begin{cor}
	
	\label{2concentrationmulitdim}
	
	Let $P,Q\in\PP_d$, $(U_t^N)_{t\in\R^+}$, $(V_t^N)_{t\in\R^+}$, $(W_t^N)_{t\in\R^+}$ be independent families of $p$ unitary Brownian motions of size $N$. Let $Z^{NM}$ be a family of deterministic matrices. With 
	$$ h : x\otimes y \in (\M_N(\C) \otimes \M_M(\C))^{\otimes 2} \mapsto y\sharp x \in \M_N(\C) \otimes \M_M(\C) ,$$
	and notations as in subsection \ref{2notationbasic}, one has for any $T\geq 0$,
	\begin{align*}
	&\E\Big[ (\tr_{N}\otimes I_{M})^{\bigotimes 2}\left( \widetilde{P}(U_T^N\otimes I_M,Z^{NM})\otimes \widetilde{Q}(U_T^N\otimes I_M,Z^{NM}) \right) \Big] \\
	&-  \E[\tr_{N}\otimes I_{M}]^{\bigotimes 2}\left( \widetilde{P}(U_T^N\otimes I_M,Z^{NM})\otimes \widetilde{Q}(U_T^N\otimes I_M,Z^{NM}) \right) \\
	&= - \frac{1}{N} \sum_{1\leq k\leq p} \int_{0}^{T} \E\Big[ \tr_N\otimes I_M\Big( h\circ \delta_k\widetilde{P}(V_{T-t}^N U_t^N\otimes I_M, Z^{NM}) h\circ \delta_k\widetilde{Q}(W_{T-t}^N U_t^N\otimes I_M, Z^{NM}) \Big) \Big] dt .
	\end{align*}
	Besides if $Z^{NM} = (I_N\otimes Y_1^M,\dots ,I_N\otimes Y_q^M)$ and that these matrices commute, then we have the same equality but with $\D_k$ instead of $h\circ \delta_k$.
	
\end{cor}

\begin{proof}
	
	Let $A^N$ be a family of deterministic matrices, by polarization and the fact that $ \D_k(Q^*) ^* = - \D_k Q $, we have
	\begin{align*}
	&\E\Big[\tr_N\left(\widetilde{P}(U_T^N,A^N)\right) \tr_N\left(\widetilde{Q}(U_T^N,A^N)\right)\Big] - \E\Big[\tr_N\left(\widetilde{P}(U_T^N,A^N)\right)\Big] \E\Big[\tr_N\left(\widetilde{Q}(U_T^N,A^N)\right)\Big] \\
	&= \E\Bigg[ \Big( \tr_N\left(\widetilde{P}(U_T^N,A^N)\right) - \E\Big[\tr_N\left(\widetilde{P}(U_T^N,A^N)\right)\Big] \Big) \\
	&\quad\quad \times \overline{\Big( \tr_N\left(\widetilde{Q^*}(U_T^N,A^N)\right) -  \E\Big[\tr_N\left(\widetilde{Q^*}(U_T^N,A^N)\right)\Big] \Big)} \Bigg] \\
	&= \frac{1}{N} \sum_{k\leq p} \int_{0}^{T} \E\Big[ \tr_N\Big( \widetilde{\D_k P}(V_{T-t}^N U_t^N, A^N) \times \widetilde{\D_k Q^*}(W_{T-t}^N U_t^N, A^N) ^* \Big) \Big] dt \\
	&= - \frac{1}{N} \sum_{k\leq p} \int_{0}^{T} \E\Big[ \tr_N\Big( m\circ \delta_k \widetilde{P}(V_{T-t}^N U_t^N, A^N) \times m\circ \delta_k \widetilde{Q}(W_{T-t}^N U_t^N, A^N) \Big) \Big] dt .
	\end{align*}
	
	\noindent Now we want to study a polynomial in $(U_T^N\otimes I_M,Z^{NM})$ and their adjoints. By linearity we can assume that $P$ is a monomial. One can also assume that $Z_i^{NM}$ is a simple tensor, i.e. that $Z_i^{NM} = A_i\otimes B_i$ where $A_i\in\M_N(\C)$ and $B_i\in\M_M(\C)$. Indeed if $Z_i^{NM} = \sum_l A_{i,l}\otimes B_{i,l}$, then a polynomial in $(U_T^N\otimes I_M,Z^{NM})$ is a linear combination in monomials in $(U_T^N\otimes I_M,(A_{i,l}\otimes B_{i,l})_{i,l})$. Thus
	\begin{align*}
	\widetilde{P}(U_T^N\otimes I_M,Z^{NM}) = \widetilde{P}(U_T^N,A) \otimes  \widetilde{P}(I_M,B) .
	\end{align*}
	
	\noindent Thus assuming that $P$ and $Q$ are monomials, we have
	
	\begin{align*}
	&\E\Big[ (\tr_{N}\otimes I_{M})^{\bigotimes 2}\left( \widetilde{P}(U_T^N\otimes I_M,Z^{NM})\otimes \widetilde{Q}(U_T^N\otimes I_M,Z^{NM}) \right) \Big] \\
	&-  \E[\tr_{N}\otimes I_{M}]^{\bigotimes 2}\left( \widetilde{P}(U_T^N\otimes I_M,Z^{NM})\otimes \widetilde{Q}(U_T^N\otimes I_M,Z^{NM}) \right) \\
	&= \Big(\E[\tr_N(\widetilde{P}(U_T^N,A)) \tr_N(\widetilde{Q}(U_T^N,A))] - \E[\tr_N(\widetilde{P}(U_T^N,A))]\ \E[ \tr_N(\widetilde{Q}(U_T^N,A))] \Big) \\
	&\quad\quad \otimes \widetilde{P}(I_M,B) \widetilde{Q}(I_M,B)\\
	&= - \frac{1}{N} \int_0^T \E\Big[ \tr_N\Big( m\circ \delta_k \widetilde{P}(V_{T-t}^N U_t^N ,A) \times m\circ \delta_k \widetilde{Q}(W_{T-t}^N U_t^N,A) \Big) \Big] \otimes \widetilde{P}(I_M,B) \widetilde{Q}(I_M,B) \ dt\\
	&= -\frac{1}{N} \sum_{k\leq p} \int_{0}^{T} \E\Big[ \tr_N\otimes I_M\Big( h\circ \delta_k\widetilde{P}(V_{T-t}^N U_t^N\otimes I_M, Z^{NM}) \times h\circ \delta_k\widetilde{Q}(W_{T-t}^N U_t^N\otimes I_M, Z^{NM}) \Big) \Big] dt .
	\end{align*}
	
	\noindent Hence the conclusion by linearity. Besides if $Z^{NM} = (I_N\otimes Y_1^M,\dots ,I_N\otimes Y_q^M)$, and that these matrices commute, then for any $p$-tuple of unitary matrices $U$, with $Z = (Y_i^M)_i$, 
	\begin{align*}
	h\circ \delta_k\widetilde{P}(U\otimes I_M, Z^{NM}) =& \sum_{P=SY_kT} \widetilde{T}(U,I_N)\widetilde{S}(U,I_N)U_k\otimes \widetilde{S}(I_M,Z)\widetilde{T}(I_M,Z) \\
	&- \sum_{P=SY_{k+d}T} U_k^*\widetilde{T}(U,I_N)\widetilde{S}(U,I_N)\otimes \widetilde{S}(I_M,Z)\widetilde{T}(I_M,Z) \\
	=& \sum_{P=SY_kT} (\widetilde{T}(U,I_N)\otimes \widetilde{T}(I_M,Z))(\widetilde{S}(U,I_N)\otimes \widetilde{S}(I_M,Z)) \times U_k\otimes I_M \\
	&- \sum_{P=SY_{k+d}T} U_k^*\otimes I_M\times (\widetilde{T}(U,I_N)\otimes \widetilde{T}(I_M,Z))(\widetilde{S}(U,I_N)\otimes \widetilde{S}(I_M,Z)) \\
	=&\ \D_k\widetilde{P}(U\otimes I_M, Z^{NM})
	\end{align*}
	Hence once again the conclusion by linearity.
	
\end{proof}

\subsection{Convergence of the free unitary Brownian motion}

If $u_t$ is a free unitary Brownian motion at time $t$, one can define $\mu_{u_t}$ as in Definition \ref{2freeprob}. Then thanks to Riesz theorem, there is a measure $\nu_t$ such that for any polynomial $P$ in two commuting variables,
$$ \tau(P(u_t,u_t^*)) = \int_{\C} P(z,z^*)\ d\nu_t(z) .$$

\noindent The measure $\nu_t$ is well-known albeit not explicit. The proof of the following theorem can be found in \cite{freebrlaw}.

\begin{theorem}
	For every $t>0$, the measure $\nu_t$ is absolutely continuous with respect to the Haar measure on $\mathbb{T} = \{z\in\C\ |\ |z|=1\}$. For $t>4$, the support of $\nu_t$ is equal to $\mathbb{T}$, and its density is positive on $\mathbb{T}$. We set $\kappa(t,\omega)$ the density of $\nu_t$ at the point $\omega\in\mathbb{T}$. Then for $t>4$, $\kappa(t,\omega)$ is the real part of the only solution with positive real part of the equation,
	\begin{equation}
	\label{2horreq}
	\frac{z-1}{z+1} e^{\frac{t}{2}z}	= \omega .
	\end{equation}
	
\end{theorem}

The following theorem states that given a free unitary Brownian motion $(u_t)_{t\geq 0}$, there exists a family of Haar unitaries $(\tilde{u}_t)_{t\geq 0}$ such that when $t$ goes to infinity, $u_t$ and $\tilde{u}_t$ are exponentially close for the operator norm topology. And more importantly it gives explicit estimates.

\begin{prop}
	\label{2Browniantounitary}
	There exists a $\CC^*$-algebra $\CC$ which contains $u_t^1,\dots,u_t^p$ freely independent free unitary Brownian motions at time $t\geq5$ and $\tilde{u}^1_t,\dots,\tilde{u}^p_t$ freely independent Haar unitaries such that for any $i$, $\norm{u_t^i - \tilde{u}_t^i} \leq 4e^2\pi e^{-\frac{t}{2}}$.
\end{prop}

\begin{proof}
	We view $B(L^2([0,2\pi]))$ as the $\CC^*$-algebra endowed with the state 
	$$\tau(u) = \langle u(\1_{[0,2\pi]}), \1_{[0,2\pi]}\rangle_{L^2([0,2\pi]}.$$
	The endomorphism $x: f\mapsto (t\to tf(t))$ is self-adjoint and has distribution (as defined in \ref{2freeprob}) $\mu_x(f) = \int_{[0,2\pi]} f$. Consequently we set $g:s\to \kappa(t,e^{is})$ and $G:s\to \int_{0}^s g(u)\ du$. Since $g$ is positive, we can define $u_t = e^{\i G^{-1}(x)}$ which has the distribution of a free unitary Brownian motion at time $t$, indeed for any polynomial $P$ in two commuting indeterminates,
	$$ \tau(P(u_t,u_t^*)) = \int_0^{2\pi} P\left(e^{\i G^{-1}(s)},e^{-\i G^{-1}(s)}\right) ds = \int_0^{2\pi} P(e^{\i s},e^{-\i s})) g(s) ds = \int_{\C} f(z,z^*)\ d\nu_t(z) .$$
	
	\noindent And similarly, $u=e^{\i x}$ is a Haar unitary. Besides, since
	$$ u_t - u = \int_0^1 e^{\i \alpha G^{-1}(x)} (G^{-1}(x)-x) e^{\i (1-\alpha)x} d\alpha ,$$
	
	\noindent thanks to the fact that $G$ is a diffeomorphism of $[0,2\pi]$,
	$$ \norm{u_t - u} \leq \norm{G^{-1}(x)-x} = \sup_{s\in [0,2\pi]} |G^{-1}(s)-s| = \sup_{s\in [0,2\pi]} |s-G(s)| \leq 2\pi \sup_{s\in [0,2\pi]} |1-g(s)| .$$
	
	\noindent We set $y(s)$ the imaginary part of the only solution with positive real part of the equation \eqref{2horreq}. Then we have for any $s$,
	$$ \frac{(g(s)-1)^2 + y(s)^2}{(g(s)+1)^2 + y(s)^2} \leq  e^{-t g(s)} .$$
	
	\noindent However since $(g(s)-1)^2\leq (g(s)+1)^2$, we have,
	$$ \frac{(g(s)-1)^2}{(g(s)+1)^2}\leq \frac{(g(s)-1)^2 + y(s)^2}{(g(s)+1)^2 + y(s)^2} \leq  e^{-t g(s)} .$$
	
	\noindent First in the case where $g(s)\geq 1$, then since we assumed $t\geq 4$, $|g(s)-1| \leq (|g(s)-1|+2) e^{-2 |g(s)-1|} e^{-\frac{t}{2}} $, and since the function $u\to (u+2)e^{-2u}$ is decreasing, we have, $|g(s)-1| \leq 2 e^{-\frac{t}{2}} $. 
	
	\noindent If $g(s)\leq 1$, then after studying the graph of the function $h:g\mapsto e^{-tg/2} -\frac{1-g}{1+g}$, we have that this function takes value $0$ in in $0$, is negative on $(0,c_{t})$ for some $c_{t}\in (0,1)$, and finally is positive for $g> c_{t}$. Since we know that $g(s)$ is positive for $t>4$ and $h(g(s))\geq 0$, if we find $g$ such that $h(g)\leq0$, then $g(s)\geq g$. Besides for $t\geq 5$, we have that $h\left(\ln(t/2)\frac{2}{t}\right)  \leq 0$. Thus necessarily $g(s) \geq \ln(t/2)\frac{2}{t}$, consequently since $g(s)\leq 1$, we know that $1- g(s) \leq 2 e^{-\frac{t}{2}g(s)}$. Hence,
	$$ 1- g(s) \leq 2 e^{-\frac{t}{2}\times \ln(t/2)\frac{2}{t}} = \frac{4}{t} . $$
	
	\noindent Thus by bootstrapping, for any s,
	$$ 1- g(s) \leq 2 e^{-\frac{t}{2}\left(1 -\frac{4}{t}\right)} = 2 e^2 e^{-\frac{t}{2}} .$$	
	
	\noindent Consequently $ \norm{u_t - u} \leq 4 e^2\pi e^{-t/2}$, and thanks to Theorem \ref{2freesum}, to conclude we just need to take the free product of $p$ copies of $B(L^2([0,2\pi]))$.
	
\end{proof}

\subsection{Free stochastic calculus and free unitary Brownian motion}

In this subsection, we consider $u_t^N = (U_1^N u_t^1,\dots ,U_p^N u_t^p)$, i.e. free unitary Brownian started in $U^N$. As we will see later, thanks to Proposition \ref{2Browniantounitary}, this will let us interpolate between $U^N = (U_1^N,\dots ,U_p^N)$ random Haar unitary matrices and $u = (u^1,\dots ,u^p)$ free Haar unitaries. Concretely if $P\in\PP_d$, we set
$$ H(t) = \tau_{N}\otimes\tau_{M}\left(\widetilde{P}(u_t^N\otimes I_M,Z^{NM})\right) .$$

\noindent Then, 
$$ H(0) = \frac{1}{NM}\tr_{MN}\left(\widetilde{P}(U^N\otimes I_M,Z^{NM})\right) , $$
$$ H(\infty) = \lim_{t\to\infty} H(t) = \tau_{N}\otimes\tau_{M}\left(\widetilde{P}(u\otimes I_M,Z^{NM})\right) . $$

\noindent To prove the second line, that is to prove that $H$ converges towards $H(\infty)$, we first use Proposition \ref{2Browniantounitary} to prove that $H$ converges towards $\tau_{N}\otimes\tau_{M}\left(\widetilde{P}((U^Nu)\otimes I_M,Z^{NM})\right)$. Then thanks to the invariance of the distribution of a Haar unitary by multiplication by a unitary operator, almost surely this quantity is equal to $ H(\infty)$. The invariance can easily be proved by using the fact (see Theorem 5.4.10 of \cite{alice}) that if $V^{kN} = (V_1^{kN},\dots,V_p^{kN})$ are independent Haar unitary matrices of size $kN$, then for any polynomial $P$,
$$ \tau(\widetilde{P}(U^Nu)) = \lim_{k\to \infty} \E\left[ \tau_{N}\otimes\tau_k(\widetilde{P}(U^N\otimes I_k \times V^{kN})) \right] = \lim_{k\to \infty} \E\left[ \tau_{kN}(\widetilde{P}(V^{kN})) \right] = \tau(\widetilde{P}(u)) .$$

\noindent Consequently as long as the integral is well-defined, we can write
$$ \frac{1}{NM}\tr_{MN}\left(\widetilde{P}(U^N\otimes I_M,Z^{NM})\right) - \tau_{N}\otimes\tau_{M}\left(\widetilde{P}(u\otimes I_M,Z^{NM})\right) = \int_0^{\infty} \frac{dH}{dt}(t)\ dt. $$

\noindent Hence we need to compute the differential of $H$ with respect to $t$, which we do in the following proposition.

\begin{prop}
	\label{2differential}
	Let the following objects be given,
	\begin{itemize}
		\item $u=(u_t^1,\dots ,u_t^p)_{t\geq 0}$ a family of $p$ free unitary Brownian motions,
		\item $U^N = (U_1^N,\dots,U_p^N)$ matrices of size $N$,
		\item $u_t^N = (U_1^N u_t^1,\dots ,U_p^N u_t^p)$ elements of $\A_N$,
		\item $Z^{NM} = (Z^{NM}_{p+1},\dots,Z^{NM}_d) $ matrices in $\M_N(\C)\otimes \M_M(\C)$,
		\item $P\in\PP_{d}$.
	\end{itemize}
	With notation as in subsection \ref{2notationbasic}, the map $H: t \mapsto \tau_N\otimes\tau_M \left(\widetilde{P}\left(u_t^N\otimes I_M,Z^{NM}\right)\right)$ is differentiable on $\R^+$ and,
	\begin{equation*}
	\frac{dH}{dt}(t) = - \frac{1}{2} \sum_{1\leq i\leq p} \tau_M\left( (\tau_N\otimes I_M) \bigotimes (\tau_N\otimes I_M)\left(\delta_i \D_i \widetilde{P}\left(u_t^N\otimes I_M,Z^{NM}\right)\right)\right).
	\end{equation*}
	
\end{prop}

\begin{proof}
	We want to use Theorem \ref{2ito} to write $H$ as an integral which we can then easily differentiate. We need to define $X_0\in\A^d$, $K$ such that for any $t\geq 0$, $K\in (L^1([0,t],\A))^d$, $U$ such that for any $t\geq 0$, $(\1_{s\leq t} U^i_s)_{t\in\R^+} \in (\B_a^{\infty})^d$, and then,
	$$ X_t = X_0 + \int_{0}^{t} K_s ds + \sum_i \int_{0}^{t} U^i_s \# dS^i_s .$$
	
	\noindent By using the linearity of the trace and the non-commutative differential, we can assume that $Z_i^{NM} = A_i\otimes B_i$ where $A_i\in\M_N(\C)$ and $B_i\in\M_M(\C)$. We then set $X_t = (u_t^N,{u_t^N}^*,A,A^*)$. Since $(A,A^*)$ is free from $\A$, the processes $K$ and $U$ associated to $(A,A^*)$ are zero. As for $(u_t^N,{u_t^N}^*)$, by definition of a free unitary Brownian motion, we have 
	$$ \forall t\geq 0, \quad  u_t^N = U^N - \int_0^t \frac{u_s^N}{2}\ ds + \i \int_0^t (u_s^N\otimes 1_{\A})\ \# dS_s , $$
	$$ \forall t\geq 0, \quad  (u_t^N)^* = (U^N)^* - \int_0^t \frac{(u_s^N)^*}{2}\ ds - \i \int_0^t (1_{\A} \otimes (u_s^N)^*)\ \# dS_s . $$
	
	\noindent To minimize cumbersome notations, for the rest of this proof we will forget about the $N$ in $u_t^N$, and assimilate $u_t^N$ with $u_t$. Consequently we set for any $s\geq 0$,
	
	\begin{align*}
	&\forall i\in [1,p], \forall j\in [1,p], \quad K_s^j = -u_{j,s}/2,\quad  U^{i,j}_s = \i\ \1_{i=j} u_{j,s}\otimes 1_{\A}, \\
	&\forall i\in [1,p], \forall j\in [p+1,2p], \quad K_s^j = -u_{j,s}^*/2, \quad U^{i,j}_s = -\i\ \1_{i=j} 1_{\A}\otimes u_{j,s}^* ,\\
	&\forall i\in [1,p], \forall j>2p, \quad K_s^j = 0, \quad U^{i,j}_s = 0\otimes 0.
	\end{align*}
	
	\noindent Thus we have for any monomial $Q$,
	
	$$ \partial Q (X) \# K\ = -\frac{1}{2}\ \sum_{1\leq i\leq p}\ \partial_i Q(X) \# u_i + \partial_i^* Q(X) \# (u_i)^* , $$
	
	\begin{align*}
	\Delta_U(Q) (X)= - \sum_{1\leq i\leq p}&\quad \langle\langle\  (\	 (\partial_i \otimes I)\circ \partial_i Q(X)\ ) \# (u_i\otimes \1_{\A},u_i\otimes \1_{\A}) \  \rangle\rangle \\
	&- \langle\langle\  (\ (\partial_i \otimes I)\circ \partial_i^* Q(X)\ ) \# (u_i\otimes \1_{\A},\1_{\A}\otimes (u_i)^* \  \rangle\rangle \\
	&- \langle\langle\  (\ (\partial_i^* \otimes I)\circ \partial_i Q(X)\ ) \# (\1_{\A}\otimes (u_i)^*,u_i\otimes \1_{\A}) \  \rangle\rangle \\
	&+ \langle\langle\  (\ (\partial_i^* \otimes I)\circ \partial_i^* Q(X)\ ) \# (\1_{\A}\otimes (u_i)^*,\1_{\A}\otimes (u_i)^*) \  \rangle\rangle . \\
	\end{align*}
	
	\noindent And thanks to Theorem \ref{2ito}, we have for any $t\geq 0$,
	$$ Q(X_t) = Q(X_0) + \int_0^t \partial Q (X_s) \# K_s\  ds + \sum_{1\leq i\leq p} \int_0^t (\partial Q(X) \sharp U^i_s)\  \#dS_s^i + \int_{0}^t \Delta_U(Q) (X_s)\  ds .$$
	
	\noindent Thus if we fix $t\in\R^+$, then for any $\varepsilon \geq -t$,
	$$ Q(X_{t+\varepsilon}) - Q(X_t) = \int_t^{t+\varepsilon} \partial Q (X_s) \# K_s\  ds + \sum_{1\leq i\leq p} \int_t^{t+\varepsilon} ( \partial Q(X) \sharp U^i_s )\  \#dS_s^i + \int_t^{t+\varepsilon} \Delta_U(Q) (X_s)\  ds .$$
	
	\noindent As we said in section \ref{2mart}, $\left( \sum_i \int_0^t (\partial Q(X) \sharp U^i_s)\  \#dS_s^i \right)_{t\geq 0}$ is a martingale, thus
	$$ \tau_N(Q(X_{t+\varepsilon})) - \tau_N(Q(X_t)) = \int_t^{t+\varepsilon} \tau_N(\partial Q (X_s) \# K_s)\  ds + \int_t^{t+\varepsilon} \tau_N(\Delta_U(Q) (X_s))\  ds .$$
	
	\noindent Finally we have,
	\begin{equation}
	\label{2derivinter}
	\frac{d \tau_N(Q(X_t))}{dt} = \tau_N(\partial Q (X_t) \# K_t) + \tau_N(\Delta_U(Q) (X_t)) .
	\end{equation}
	
	\noindent Besides,
	$$ \tau_N(\partial Q (X_t) \# K_t) = -\frac{1}{2} \sum_{1\leq i\leq p}\ \tau_N(D_i Q(X_t)\ u_{i,t}) + \tau_N(u_{i,t}^* D_i^*P(X_t)) , $$
	
	\noindent and,
	$$ \tau_N\left(\langle\langle\  (\partial_i \otimes I)\circ \partial_i Q(X) \# (u_i\otimes \1_{\A},u_i\otimes \1_{\A}) \  \rangle\rangle\right) = \sum_{Q = AY_iBY_iC} \tau_N(C(X)A(X) u_i)\tau_N(u_i B(X)) ,$$
	$$ \tau_N\left( \langle\langle\  (\partial_i \otimes I)\circ \partial_i^* Q(X) \# (u_i\otimes \1_{\A},\1_{\A}\otimes u_i^* \  \rangle\rangle \right) = \sum_{Q = AY_iBY_i^*C} \tau_N(A(X) u_i u_i^* C(X))\tau_N(B(X)) , $$
	$$ \tau_N\left( \langle\langle\  (\partial_i^* \otimes I)\circ \partial_i Q(X) \# (\1_{\A}\otimes u_i^*,u_i\otimes \1_{\A}) \  \rangle\rangle \right) = \sum_{Q = AY_i^*BY_iC} \tau_N(A(X)C(X))\tau_N(u_i u_i^* B(X)) , $$
	$$\tau_N\left( \langle\langle\  (\partial_i^* \otimes I)\circ \partial_i^* Q(X) \# (\1_{\A}\otimes u_i^*,\1_{\A}\otimes u_i^*) \  \rangle\rangle \right) = \sum_{Q = AY_i^*BY_i^*C} \tau_N(C(X)A(X)u_i^*)\tau_N(u_i^*B(X)) .$$
	
	\noindent Besides we also have,
	
	\begin{align*}
	\tau_N\otimes\tau_N\left( \delta_i D_i Q(X) \times 1\otimes u_i \right) =&\ 2 \sum_{Q = AY_iBY_iC} \tau(B(X)u_i)\tau(C(X)A(X)u_i) \\
	&- \sum_{Q = AY_i^*BY_iC} \tau(C(X)A(X))\tau(B(X)u_i u_i^*) \\
	&- \sum_{Q = AY_iBY_i^*C} \tau(B(X))\tau(C(X)u_i u_i^*A(X)) ,
	\end{align*}
	\begin{align*}
	\tau_N\otimes\tau_N\left( \delta_i D_i^* Q(X) \times u_i^*\otimes 1 \right) =& -2 \sum_{Q = AY_i^*BY_i^*C} \tau(u_i^* B(X))\tau(C(X)A(X) u_i^*) \\
	&+ \sum_{Q = AY_i^*BY_iC} \tau(B(X)u_i u_i^*)\tau(C(X)A(X)) \\
	&+ \sum_{Q = AY_iBY_i^*C} \tau(C(X)u_i u_i^*A(X))\tau(B(X)) .
	\end{align*}
	
	\noindent Which means that
	$$ \tau_N\left( \Delta_U(Q) (X) \right) = -\frac{1}{2} \sum_{1\leq i\leq p} \tau_N\otimes\tau_N\left( \delta_i D_i Q(X) \times 1\otimes u_i \right) - \tau_N\otimes\tau_N\left( \delta_i D_i^* Q(X) \times u_i^*\otimes 1 \right) .$$
	
	\noindent And thus when combined with equation \eqref{2derivinter}, we get that
	\begin{align*}
	\frac{d \tau_N(Q(X_t))}{dt} = -\frac{1}{2}\ \sum_{1\leq i\leq p}\ &\tau_N(D_i Q(X_t)\ u_{i,t}) + \tau_N\otimes\tau_N\left( \delta_i D_i Q(X_t) \times 1\otimes u_{i,t} \right) \\
	&  + \tau_N(u_{i,t}^* D_i^*Q(X_t)) - \tau_N\otimes\tau_N\left( \delta_i D_i^* Q(X_t) \times u_{i,t}^*\otimes 1 \right) \\
	= -\frac{1}{2}\ \sum_{1\leq i\leq p}\ &\tau_N\otimes\tau_N\left( \delta_i \left(D_i Q(X_t)u_{i,t}\right) \right) - \tau_N\otimes\tau_N\left( \delta_i \left(u_{i,t}^* D_i^* Q(X_t)\right) \right) \\
	= -\frac{1}{2}\ \sum_{1\leq i\leq p}\ &\tau_N\otimes\tau_N\left( \delta_i \D_i Q (X_t) \right) .
	\end{align*}
	
	\noindent Now we want to study a polynomial in $(u_t^N,Z^{NM})$ and their adjoints. If $P$ is a monomial, we have,
	$$ \widetilde{P}(u_t^N\otimes I_M,Z^{NM}) = \widetilde{P}(u_t^N,A)\otimes \widetilde{P}(I_M,B) .$$
	
	\noindent Therefore,
	\begin{align*}
	&\frac{dH}{dt}(t) = -\frac{1}{2} \sum_{1\leq i\leq p} \tau_N\otimes\tau_N\left(\delta_i \D_i \widetilde{P}(u_t^N,A)\right) \times \tau_M\left( \widetilde{P}(I_M,B) \right) .
	\end{align*}
	
	\noindent And since for any $S,T\in \PP_d$, 
	\begin{align*}
	&\tau_N\otimes\tau_N(\delta_i \widetilde{TS}(u_t^N,A)) \times \tau_M\left( \widetilde{ST}(I_M,B) \right) \\
	&= \tau_{M}\left( (\tau_N\otimes I_M)\bigotimes (\tau_N\otimes I_M)\left( \delta_i \widetilde{TS}\left(u_t^{N}\otimes I_M, A\otimes B\right) \right) \right) .
	\end{align*}
	
	\noindent Hence after summing,
	\begin{align*}
	&\frac{d}{dt}\tau_N\otimes\tau_M \left(\widetilde{P}\left(u_t^N,Z^{NM}\right)\right) = -\frac{1}{2} \sum_{1\leq i\leq p} \tau_M\left( (\tau_N\otimes I_M) \bigotimes (\tau_N\otimes I_M)\left(\delta_i \D_i \widetilde{P}\left(u_t^N,Z^{NM}\right)\right)\right) ,
	\end{align*}
	and we conclude by linearity.	
	
\end{proof}

\section{Proof of Theorem \ref{2imp1}, the main result}
\label{2mainsec}

\subsection{Overview of the proof}
\label{2overview}

If we take the point of view of free probability -- for details we refer to the third point of Definition \ref{2freeprob} -- we have two families of non-commutative random variables, $(U^N\otimes I_M, Z^{NM})$ and $(u\otimes I_M, Z^{NM})$, and we want to study the difference between their distributions. As mentioned in the introduction the main idea of the proof is to interpolate those two families with the help of $p$ free unitary Brownian motions $(u^N_t)_{t\geq 0}$ started in the Haar unitary matrices $U^N$. A big difference with the case of the GUE which was treated in \cite{un} is that we do not have an explicit expression of the law of $u^N_t$ in function of $U^N$ and $u$, which is why we had to introduce notions of free stochastic calculus.

This idea of interpolating random matrices is not entirely new. Indeed in \cite{UBMN}, the authors proved in theorem 1.3 that given $(U^N_t)_{t\geq 0}$ a unitary Brownian motion of size $N$ and $(u_t)_{t\geq 0}$ a free unitary Brownian motion, for any $n\in\N$,
\begin{equation}
\label{2nvrapport}
\left| \E\left[\tau_{N}((U^N_t)^n)\right] - \tau(u_t^n) \right| \leq \frac{t^2n^4}{N^2}.
\end{equation}
The main idea of their proof was to interpolate $\E\left[\tau_{N}((U^N_t)^n)\right]$ and $\E\left[\tau_{N}((U^{2N}_t)^n)\right]$ through a stochastic process defined with a unitary Brownian motion. Thus they get an estimate of the difference between those two expectations, and by iterating this method they get an estimate of 
$$\E\left[\tau_{N}((U^{2^iN}_t)^n)\right] - \E\left[\tau_{N}((U^{2^{i+1}N}_t)^n)\right]$$
for any $i$. Then by summing over $i$ they get equation \eqref{2nvrapport}. In this paper, while the method used are very different, we do keep this idea of interpolation. But instead of doing it step by step from $2^iN$ to $2^{i+1}N$, we directly interpolate between $N$ and $\infty$.

Since our aim in this subsection is not to give a proof but to outline the strategy used in subsection \ref{2technical}, we assume that we have no matrix $Z^{NM}$ and that $M=1$. Now under the assumption that this is well-defined, if $Q$ is a non-commutative polynomial,

$$	\E\left[\frac{1}{N}\tr_{N}\Big(Q\left(U^N\right)\Big)\right] - \tau\Big(Q\left(u\right)\Big) = - \int_{0}^{\infty} \E\left[\frac{d}{dt}\Big(\tau_N\left(Q(u^N_t)\right)\Big)\right]\ dt .$$

\noindent In the classical case, if $(S_t)_{t\geq 0}$ is a Markov process with generator $\theta$, then under the appropriate assumption we have
$$ \frac{d}{dt} \E[f(S_t)] = \E[(\theta f)(S_t)] . $$

\noindent And if the law of the process at time $0$ is invariant for this Markov process we have that for any $t$, $\E[(\theta f)(S_t)]=0$. Since $(u_t^N)_{t\geq 0}$ is a free Markov process, we expect to get similarly that

$$ \frac{d}{dt}\Big(\tau_N\left(Q(u^N_t)\right)\Big) = \tau_N\left((\Theta Q)(u_t^N)\right) , $$

\noindent for some generator $\Theta$ which we will compute with the help of Proposition \ref{2differential}. Besides the invariant law of a free Brownian motion is the law of free Haar unitaries. Thus if $(u_t)_{t\geq 0}$ is a free Brownian motion started in free Haar unitaries, we have that $\tau\left((\Theta Q)(u_t)\right) = 0$. Since unitary Haar matrices converges in distribution towards free Haar unitaries (see \cite{alice}, Theorem 5.4.10), we have that $\tau_N\left((\Theta Q)(u_t^N)\right)$ converges towards $\tau\left((\Theta Q)(u_t)\right) = 0$. As we will see in this proof, the convergence happens at a speed of $N^{-2}$. To prove this, the main idea is to view free unitary Brownian motions started in $U^N$ as the asymptotic limit when $k$ goes to infinity of a unitary Brownian motion of size $kN$ started in $U^N\otimes I_k$ (see Proposition \ref{2limite_trace}).

Another issue is that to prove Theorem \ref{2imp1}, we would like to set $Q = f(P)$ but since $f$ is not polynomial this means that we need to extend the definition of operators such as $\delta_i$. In order to do so we assume that there exists $\mu$ a measure on $\R$ such that,
$$ \forall x\in\R,\quad f(x) = \int_{\R} e^{\i xy}\ d\mu(y) .$$

\noindent While we have to assume that the support of $\mu$ is indeed on the real line, $\mu$ can be a complex measure. However we will usually work with measure such that $|\mu|(\R)$ is finite. Indeed under this assumption we can use Fubini's theorem, and we get

$$	\E\left[\frac{1}{M}\tr_{N}\Big(f\left(P(U^N)\right)\Big)\right] - \tau\Big(f\left(P(u)\right)\Big) = \int_{\R}	\left\{ \E\left[\frac{1}{N}\tr_{N}\Big(e^{\i y P\left(U^N\right)}\Big)\right] - \tau\Big(e^{\i y P\left(u\right)}\Big)\right\}\ d\mu(y) . $$

\noindent We can then set $Q = e^{\i y P}$. And even though this is not polynomial, since it is a power series, most of the properties associated to polynomials remain true with some assumption on the convergence. The main difficulty with this method is that we need to find a bound uniform in $y$, indeed we have terms of the form
$$ \int_{\R}|y|^l\ d|\mu|(y) $$
which appear. Thanks to Fourier integration we can relate the exponent $l$ to the regularity of the function $f$, thus we want to find a bound with $l$ as small as possible. It turns out that with our proof $l = 4$.

\subsection{Proof of Theorem \ref{2imp2}}

\label{2technical}

In this section we focus on proving theorem \ref{2imp2} from which we deduce all of the important corollaries, and notably Theorem \ref{2imp1}. Since this subsection is dedicated to proving only one theorem but is the longest of the paper, in the next few paragraphs we explain the structure of the proof. Unlike subsection \ref{2overview} where we gave ideas on the method of the proof, here me mainly focus on which specific purpose serves every lemma. We begin by applying Proposition \ref{2differential} to directly obtain Lemma \ref{2fourier} which states that the difference between the trace at time $N$, that is a trace in Haar unitary matrices, and at infinity, that is a trace in free Haar unitaries, can be written as an integral with respect to $t$ of a trace in free unitary Brownian motions at time $t$. We then proceed to study the term under the integral that we name $S_{t,y}^N(A,B)$ in Definition \ref{2quantpert}. 

The next two lemmas are technical lemmas that we need to justify further computations. First Proposition \ref{2limite_trace} shows that one can see a trace of a power series in free unitary Brownian motions and matrices $U^N$ of size $N$ as the limit of the trace of the same power series but evaluated in independent unitary Brownian motions of size $kN$ and $U^N\otimes I_k$. The proof of Proposition \ref{2limite_trace} can be summarized as using well-known theorems on the convergence in distribution of random matrices. The proof of the second one, Lemma \ref{2nondiag}, is much more subtle. It gives an estimate in $k$ of the non diagonal coefficients of our power series in independent unitary Brownian motions of size $kN$ and $U^N\otimes I_k$. The proof relies on Gronwall's inequality to reduce the problem to the polynomial case and classical stochastic calculus to deal with the former.

Lemma \ref{2grcalc} let us write $S_{t,y}^N(A,B)$ as a linear combination of covariance terms. Indeed, thanks to Proposition \ref{2limite_trace} one can write $S_{t,y}^N(A,B)$ as the limit when $k$ goes to infinity of a linear combination of products of expectations of traces in power series in independent unitary Brownian motions of size $kN$ started in $U^N\otimes I_k$ where $U^N$ are Haar unitary matrices. It turns out that the expectation of the product of those traces converges towards $0$ when $k$ goes to infinity thanks to Lemma \ref{2nondiag} and the properties of Haar unitary matrices. Hence every product of expectations can be viewed as a covariance term. Interestingly enough, this is the only part of the proof where we use that $U^N$ are Haar unitary matrices. More precisely we use that the law of such a random matrix is invariant by multiplication by a unitary matrix. In every other part of the proof one only need to assume at most that $U^N$ are unitary matrices.

Finally, in Lemma \ref{2optilll}, we use Corollary \ref{2concentrationmulitdim} to get an upper bound independent of $k$ of those covariance terms. Thus we can let $k$ go to infinity to get an upper bound of $S_{t,y}^N(A,B)$. As usual, the covariance of renormalized traces on $\M_N(\C)$ is of order $N^{-2}$. And even though we are not exactly working with a trace on $\M_N(\C)$ here, when using Corollary \ref{2concentrationmulitdim} the upper bound that we get is of order $N^{-2}$. Which means that so is the difference between the trace at time $N$ and at infinity. Finally this upper bound immediately yields Theorem \ref{2imp2}. Then we conclude this section by a proof of Theorem \ref{2imp1} which, up to a trick to assume that they have compact support, mainly consist in checking that the functions that we consider in Theorem \ref{2imp1} satisfies the hypothesis of Theorem \ref{2imp2}.

\begin{theorem}
	\label{2imp2}
	We define
	\begin{itemize}
		\item $u=(u^1,\dots ,u^p)$ a family of $p$ free Haar unitaries,
		\item $U^N = (U_1^N,\dots,U_p^N)$ i.i.d. Haar unitary matrices of size $N$.
		\item $Z^{NM} = (Z_{p+1}^{NM},\dots,Z_d^{NM})$ deterministic matrices,
		\item $P\in \PP_{d}$ a self-adjoint polynomial,
		\item $f:\R\mapsto\R$ such that there exists a measure on the real line $\mu$ with $\int (1+y^4)\  d|\mu|(y)\ < +\infty$ and for any $x\in\R$,
		$$ f(x) = \int_{\R} e^{\i x y}\ d\mu(y) . $$
	\end{itemize}
	
	\noindent Then there exists a polynomial $L_P\in\R^+[X]$ which only depends on $P$ such that for any $N,M$,
	\begin{align*}
	\Bigg| &\E\left[\frac{1}{MN}\tr_{MN}\Big(f\left(\widetilde{P}\left(U^N\otimes I_M,Z^{NM}\right)\right)\Big)\right] - \tau_N\otimes\tau_M\Big(f\left(\widetilde{P}\left(u\otimes I_M,Z^{NM}\right)\right)\Big) \Bigg| \\
	&\leq \frac{\ln^2(N) M^2}{N^2} L_P\left(\norm{Z^{NM}}\right)\times \int_{\R} (|y|+y^4)\ d|\mu|(y) .
	\end{align*}
	
	\noindent where $\norm{Z^{NM}} = \sup\limits_{p+1\leq i\leq d} \norm{Z_i^{NM}}$.
	
\end{theorem}

Even though we do not give an explicit expression for $L_P$, it is possible to compute it rather easily by following the proof of Lemma \ref{2optilll}. In particular given a set of polynomials whose degree and coefficients are uniformly bounded, we can find a polynomial $R$ such that for any $P$ in this set and any matrices $Z^{NM}$, $L_P\left(\norm{Z^{NM}}\right) \leq R\left(\norm{Z^{NM}}\right)$. Besides, if we replace $P$ by $\alpha P$ where $\alpha\in\C$, then up to a constant independent from $\alpha$, we can bound $L_{\alpha P}$ by $(|\alpha|+|\alpha|^5) L_P$, or even $(|\alpha|+|\alpha|^4) L_P$ if one picks the first expression in the minimum. It is also wort noting that the set of function $f:\R\mapsto\C$ such that there exists a Borel complex measure on the real line $\mu$ with $\int y^4\  d|\mu|(y)\ < +\infty$ and for any $x\in\R$, $ f(x) = \int_{\R} e^{\i x y}\ d\mu(y) $, is the so-called $4^{\text{th}}$ Wiener space $W_4(\R)$. We refer to \cite{Evangelos}, section 4.3 for a brief introduction on the matter. 

The first step to prove this theorem is the following lemma, who is a direct consequence of Proposition \ref{2differential} and equation \eqref{2extension},

\begin{lemma}
	\label{2fourier}
	
	With the same notation as in Theorem \ref{2imp2}, we define
	\begin{itemize}
		\item $u=(u_t^1,\dots ,u_t^p)_{t\geq 0}$ a family of $p$ free unitary Brownian motions,
		\item $u_t^N = (U_1^N u_t^1,\dots ,U_p^N u_t^p)$ elements of $\A_N$.
	\end{itemize}
	Then with notation as in subsection \ref{2notationbasic}, almost surely
	\begin{align*}
	& \frac{1}{MN}\tr_{MN}\left(f\left(\widetilde{P}\left(U^N\otimes I_M,Z^{NM}\right)\right)\right) - \tau_N\otimes\tau_M\left(f\left(\widetilde{P}\left(u_T^N\otimes I_M,Z^{NM}\right)\right)\right) \\
	&= \frac{1}{2} \sum_{i\leq p} \int \int_0^T \tau_M\left( (\tau_N\otimes I_M) \bigotimes (\tau_N\otimes I_M)\left(\delta_i\left( \D_i\ e^{\i y \widetilde{P}} \right) \left(u_t^N\otimes I_M,Z^{NM}\right)\right)\right) dt\ d\mu(y) .
	\end{align*}
	
\end{lemma}

Since $\D_i\ e^{\i y P} = \i y\ \delta_i P\ \widetilde{\#}\ e^{\i y P}$, this prompts us to define the following quantity.

\begin{defi}
	\label{2quantpert}
	Let $A,B\in \PP_{d}$, we set 
	$$S_{t,y}^N(A,B) = \E\left[ \tau_M\left( (\tau_N\otimes I_M) \bigotimes (\tau_N\otimes I_M)\left(\delta_i\left( A\ e^{\i y P}\ B \right) \left(Z_t^N\right)\right)\right) \right] ,$$
	where $Z_t^N = \left(u_t^N\otimes I_M,Z^{NM},(u_t^N)^*\otimes I_M,(Z^{NM})^*\right)$.
\end{defi}

The following proposition justifies why the family $(U^N\otimes I_M, u_t\otimes I_M, Z^{NM})$ has in the large $k$ limit the distribution -- in the sense of Definition \ref{2freeprob} -- of the family $(U^N\otimes I_{kM}, U_t^{kN}\otimes I_M, Z^{NM}\otimes I_k)$ where $U_t^{kN}$ are independent unitary Brownian motions of size $kN$ at time $t$.

\begin{prop}
	\label{2limite_trace}
	
	If $U_t^{kN}$ are unitary Brownian motions of size $kN$ at time $t$, independent of $U^N$, we set
	\begin{equation*}
	Y_t^k = \Big((U^N\otimes I_k\ U_t^{kN})\otimes I_M, Z^{NM}\otimes I_k, (U^N\otimes I_k\ U_t^{kN})^*\otimes I_M, (Z^{NM})^*\otimes I_k\Big) .
	\end{equation*}
	
	\noindent Then if $q = A e^{\i y P} B$, we have that for any $t$, almost surely with respect to $U^N$,
	
	$$ (\tau_N\otimes I_M)\big(q(Z_t^N)\big) = \lim_{k\to\infty} \E_k\left[ (\tau_{kN}\otimes I_M)\big(q(Y_t^k)\big) \right] ,$$
	
	\noindent where $\E_k$ is the expectation with respect to $(U_t^{kN})_{t\geq 0}$.
	
\end{prop}

\begin{proof}
	
	It has been known for a long time that the unitary Brownian motion converges in distribution towards the free unitary Brownian motion, see \cite{bianebr}. However since we also have to consider deterministic matrices we will use Theorem 1.4 of \cite{UBMN}. This theorem states that if $(U_t^{kN})_{t\geq 0}$ are independent unitary Brownian motions and $D^{kN}$ is a family of deterministic matrices which converges strongly in distribution towards a family of non-commutative random variables $d$, the family $(U_t^{kN},D^{kN})$ in the non-commutative probability space $(\M_{kN}(\C),*,\E_k[\frac{1}{kN}\tr])$ converges strongly in distribution towards the family $(u_t,d)$ where $u_t$ are freely independent free unitary Brownian motions at time $t$ free from $d$. That being said, we do not use the convergence of the norm, we only need the convergence in distribution which is way easier to prove through induction and stochastic calculus. In our situation we can write for every $i$,
	
	$$ Z_i^{NM} = \sum_{1\leq r,s\leq N} E_{r,s}\otimes A_{r,s,i}^M .$$
	
	\noindent Thus if $E^N = (E_{r,s})_{1\leq r,s\leq N}$, we fix $D^{kN} = (U^N\otimes I_k, E^N\otimes I_k)$, and we can apply Theorem 1.4 of \cite{UBMN} to get that for any non-commutative polynomial $P$,
	
	$$ \lim_{k\to\infty} \E_k\left[\tau_{kN}(\widetilde{P}(U_t^{kN},D^{kN})) \right] = \tau_N\left(\widetilde{P}(u_t,U^N,E^N)\right) .$$
	
	\noindent Consequently, for any non-commutative polynomial $P$, we also have
	
	$$ \lim_{k\to\infty} \E_k\left[\tau_{kN}\otimes I_M (\widetilde{P}(U_t^{kN},D^{kN},A^M)) \right] = \tau_N\otimes I_M  \left(\widetilde{P}(u_t,U^N,E^N,A^M)\right) .$$
	
	\noindent Hence for any $P\in \PP_d$,
	$$ \lim_{k\to\infty} \E_k\left[\tau_{kN}\otimes I_M (P(Y_t^k)) \right] = \tau_N\otimes I_M \left(P(Z_t^N)\right) .$$
	
	\noindent Now since $U_t^{kN}$ are unitary matrices, we can find a polynomial $L$ such that for any $k$, $\norm{P(Y_t^k)}\leq C = L\left(\norm{U^N},\norm{Z^{NM}}\right)$. Knowing this, let $f_{\varepsilon}\in \C[X]$ be a polynomial which is $\varepsilon$-close to $x\mapsto e^{\i y x}$ on the interval $[-C,C]$. Since one can always assume that $C> \norm{P(Z_t^N)}$, we have a constant $K$ such that
	
	$$ \norm{ (\tau_N\otimes I_M)\big(q(Z_t^N)\big) - (\tau_{N}\otimes I_M)\big((Af_{\varepsilon}(P)B)(Z_t^N) } \leq K \varepsilon , $$
	$$ \norm{ (\tau_N\otimes I_M)\big(q(Y_t^k)\big) - (\tau_{N}\otimes I_M)\big((Af_{\varepsilon}(P)B)(Y_t^k) } \leq K \varepsilon .$$
	
	\noindent Thus
	\begin{align*}
	&\norm{ (\tau_N\otimes I_M)\big(q(Z_t^N)\big) - \E_k\left[ (\tau_{kN}\otimes I_M)\big(q(Y_t^k)\big) \right] } \\
	&\leq \norm{ (\tau_N\otimes I_M)\big((Af_{\varepsilon}(P)B)(Z_t^N)\big) - \E_k\left[ (\tau_{kN}\otimes I_M)\big((Af_{\varepsilon}(P)B)(Y_t^k)\big) \right] }  + 2 K \varepsilon .
	\end{align*}
	
	\noindent Consequently
	
	$$\limsup\limits_{k\to \infty} \norm{ (\tau_N\otimes I_M)\big(q(Z_t^N)\big) - \E_k\left[ (\tau_{kN}\otimes I_M)\big(q(Y_t^k)\big) \right] } \leq 2K \varepsilon .$$
	
	\noindent This completes the proof.
	
\end{proof}

The next lemma shows that the non-diagonal coefficients can actually be neglected.

\begin{lemma}
	\label{2nondiag}
	We define $Y_t^k$ as in Proposition \ref{2limite_trace}, $P_{1,2} = I_N\otimes E_{1,2}\otimes I_M$, $q = A e^{\i y P} B$, then
	$$  \lim_{k\to \infty} k^{1/2} \E_k\left[(\tr_{kN}\otimes I_M)(q(Y_t^k)P_{1,2})\right] = 0 .$$
	
\end{lemma}

\begin{proof}
	
	\noindent Let us first define for $A,B\in \PP_{d}$,
	
	$$ f_{A,B}^t(y) = \E_k\left[(\tr_{kN}\otimes I_M)\Big((\widetilde{A}\ e^{\i y \widetilde{P}}\ \widetilde{B})(U^N\otimes I_{kM}, U_t^{kN}\otimes I_M, Z^{NM}\otimes I_k)\ P_{1,2}\Big)\right] ,$$
	
	$$ d_n^t(y) = \sup\limits_{\substack{A,B\in\PP_{d} \text{ monomials, }\\ \deg(AB)\leq n \\  0\leq s\leq t}} \norm{f_{A,B}^t(y)} .$$
	
	\noindent Since given a matrix $Z\in \M_{NkM}(\C)$, we have 
	\begin{align*}
	\norm{(\tr_{kN}\otimes I_M)(Z P_{1,2})} &= \norm{(\tr_{N}\otimes I_M)(I_{NM}\otimes f_2^*\times Z\times I_{NM}\otimes f_1)} \\
	&\leq N\norm{I_{NM}\otimes f_2^*\times Z\times I_{NM}\otimes f_1} \\
	&\leq N \norm{Z}. 
	\end{align*}
	
	
	\noindent Consequently, we can find a constant $D$ such that for any $n$, $ d_n^t(y) \leq D^n$. Note that this constant $D$ can be exponentially large in $N$ or $M$, indeed it does not matter since in the end we will show that this quantity tends to $0$ when $k$ goes to infinity and the other parameters such as $N,M$ or $y$ are fixed. This implies that for $a$ small enough,
	$$g_{k,a}^t(y) = \sum_{n\geq 0} d_n^t(y) a^n ,$$
	
	\noindent is well-defined. But if we set $c_R(P)$ the coefficient associated to the monomial $R$ in $P$, we have for any $s\leq t$,
	
	$$ \left| \frac{d f_{A,B}^s(y)}{dy} \right| \leq \sum_{R \text{ monomials}} |c_R(P)|\ d_{\deg(AB)+ \deg(R)}^t(y)  .$$
	
	\noindent Thus if $\deg(AB)\leq n$, we have for any $y\geq 0$,
	
	$$ f_{A,B}^s(y) \leq f_{A,B}^s(0) + \sum_{R \text{ monomials}} |c_R(P)|\ \int_0^y d_{n + \deg(R)}^t(u)\  du . $$
	
	\noindent Thus we have for $y\geq 0$ and any $n\geq 0$,
	
	$$ a^n d_n^t(y) \leq a^n d_n^t(0) + \sum_{L \text{ monomials}} |c_R(P)| a^{-\deg(L)}\ \int_0^y d_{n + \deg(R)}^t(u) a^{n+ \deg(L)} du . $$
	
	\noindent And with $\norm{.}_{a^{-1}}$ defined as in \eqref{2normA}, we have
	$$ g_{k,a}^t(y) \leq g_{k,a}^t(0) + \norm{P}_{a^{-1}} \int_0^y g_{k,a}^t(u) du .$$
	
	\noindent Thanks to Gronwall's inequality, we have for $y\geq 0$,
	
	\begin{equation}
	\label{2gronwall}
	g_{k,a}^t(y) \leq g_{k,a}^t(0) e^{y\norm{P}_{a^{-1}}} .
	\end{equation}
	
	\noindent In order to conclude the proof, we are going to show that
	\begin{equation}
	\label{2touten0}
	g_{k,a}^t(0) = \mathcal{O}(k^{-2}).
	\end{equation}
	
	\noindent In combination with equation \eqref{2gronwall}, it will yields 
	$$ \left| k^{3/2} \E_k\left[(\tau_{kN}\otimes I_M)(q(Y_t^k)P_{1,2})\right] \right| \leq k^{1/2} a^{-\deg(AB)} g_{k,a_k}^t(y) = \mathcal{O}(k^{-3/2}) .$$
	
	\noindent Hence the conclusion. To show equation \eqref{2touten0}, first one can find deterministic matrices $L_j^{u,v} \in \M_N(\C)$ such that for every $j\in [p+1,d]$, $ Z_j^{NM} = \sum_{ 1\leq u,v \leq M} L_j^{u,v}\otimes E_{u,v}$ where $E_{u,v}\in\M_M(\C)$ is the matrix whose every coefficient is $0$ but the $(u,v)$ coefficients. Thus with $L = (L_j^{u,v})_{j,u,v}$, we proceed by defining
	\begin{equation}
	\label{2definul}
	V_{N,k}^t = \Big( U_t^{kN}, (U_t^{kN})^*, U^N\otimes I_k, (U^N)^*\otimes I_k, L\otimes I_k, L^*\otimes I_k \Big) ,
	\end{equation}
	
	$$ c_n^t = \sup_{ \substack{\deg(Q)\leq n,\ Q \text{ monomial} \\ 0\leq s\leq t}} \left| \E_k\left[\tr_{kN}(Q(V_{N,k}^s)\ P_{1,2})\right] \right| .$$
	
	\noindent where for the rest of the proof $P_{1,2} = I_N \otimes E_{1,2}$. Then for any $A\in\PP_d$ monomials,
	\begin{equation}
	\label{2necquote}
	(\tr_{kN}\otimes I_M)\Big(\widetilde{A}(U^N\otimes I_{kM}, U_t^{kN}\otimes I_M, Z^{NM}\otimes I_k)\ P_{1,2}\Big)
	\end{equation}
	is a linear combination of at most $M^{2n}$ terms of the form
	\begin{equation}
	\label{2necquote2}
	\tr_{kN}\Big(A_r(V_{N,k}^t)\ P_{1,2}\Big) \times E_{u,v},
	\end{equation}
	where $(A_r)_r$ are monomials. Consequently we have
	$$ d_n^t(0) \leq M^{2n} c_n^t .$$
	
	\noindent Thus if we set 
	$$f_k^t(a) = \sum_{n\geq 0} c_n^t a^n , $$ 
	
	\noindent we have 
	$$ g_{k,a}^t(0) \leq f_k^t(M^2 a) . $$
	
	\noindent So all we need to do is to prove that for $a$ small enough, $f_k^t(a) = \mathcal{O}(k^{-2})$. Let $Q$ be a monomial, we define $Q_t$ as the monomial $Q$ evaluated in $V^t_{N,k}$. Thanks to Proposition \ref{2unitarybr},
	
	\begin{align*}
	\frac{d}{dt}\E_k\left[\tr_{kN}\left( Q_t P_{1,2}\right)\right] = &- \frac{|Q|_B}{2}\ \E_k\left[\tr_{kN}\left( Q_t P_{1,2}\right) \right] \\
	&- \frac{1}{kN} \sum_{1\leq i\leq p,\  Q = A U_i B U_i C} \E_k\left[\tr_{kN}\left(A_t U^{kN}_{i,t} C_tP_{1,2}\right) \tr_{kN}\left(B_t U^{kN}_{i,t}\right)\right] \\
	&- \frac{1}{kN} \sum_{1\leq i\leq p,\  Q = A U_i^* B U_i^* C} \E_k\left[\tr_{kN}\left(A_t {U^{kN}_{i,t}}^* C_tP_{1,2}\right) \tr_{kN}\left(B_t {U^{kN}_{i,t}}^*\right)\right] \\
	&+ \frac{1}{kN} \sum_{1\leq i\leq p,\ Q=A U_i B U_i^* C } \E_k\left[\tr_{kN}(A_t C_tP_{1,2}) \tr_{kN}(B_t)\right] \\
	&+ \frac{1}{kN} \sum_{1\leq i\leq p,\ Q=A U_i^* B U_i C } \E_k\left[\tr_{kN}(A_t C_tP_{1,2}) \tr_{kN}(B_t)\right] .
	\end{align*}
	
	\noindent Since $\E_k\left[\tr_{kN}\left( Q_0 P_{1,2}\right)\right] = 0$, we have for any $t$,
	
	\begin{align*}
	\E_k\left[\tr_{kN}\left( Q_t P_{1,2}\right)\right]& \\
	= \int_{0}^t e^{- \frac{|Q|_B}{2} (t-s)} \Bigg( &- \frac{1}{kN} \sum_{1\leq i\leq p,\  Q = A U_i B U_i C} \E_k\left[\tr_{kN}\left(A_s U^{kN}_{i,s} C_sP_{1,2}\right) \tr_{kN}\left(B_s U^{kN}_{i,s}\right)\right] \\
	&- \frac{1}{kN} \sum_{1\leq i\leq p,\  Q = A U_i^* B U_i^* C} \E_k\left[\tr_{kN}\left(A_s {U^{kN}_{i,s}}^* C_tP_{1,2}\right) \tr_{kN}\left(B_s {U^{kN}_{i,s}}^*\right)\right] \\
	&+ \frac{1}{kN} \sum_{1\leq i\leq p,\ Q=A U_i B U_i^* C } \E_k\left[\tr_{kN}(A_s C_sP_{1,2}) \tr_{kN}(B_s)\right] \\
	&+ \frac{1}{kN} \sum_{1\leq i\leq p,\ Q=A U_i^* B U_i C } \E_k\left[\tr_{kN}(A_s C_sP_{1,2}) \tr_{kN}(B_s)\right] \Bigg) ds .
	\end{align*}
	
	\noindent As in Proposition \ref{2concentration}, we consider $(V_t^{kN})_{t\in\R^+}$ and $(W_t^{kN})_{t\in\R^+}$ independent families of $p$ unitary Brownian motions of size $kN$, independent of $(U_t^{kN})_{t\in\R^+}$. We define $V_{N,k}^{r,1}$ and $V_{N,k}^{r,2}$ as $V_{N,k}^r$ (see \eqref{2definul}) but with $V_{s-r}^{kN} U_r^{kN}$ and $W_{s-r}^{kN} U_r^{kN}$ instead of $U_r^{kN}$. Thanks to Proposition \ref{2concentration}, by polarization and the fact that $(\D_i Q^*)^* = - \D_i Q$, we have with $\cov(X,Y) = \E[XY] - \E[X]\E[Y]$,
	\begin{align*}
	&\cov_k\left( \tr_{kN}(A_s C_sP_{1,2}), \tr_{kN}(B_s) \right) \\
	&= -\frac{1}{kN} \sum_{i\leq p} \int_{0}^{s} \E_k\Big[ \tr_{kN}\Big( (\delta_i(AC) \widetilde{\#} P_{1,2})\left(V_{N,k}^{r,1}\right)\ (\D_i B)\left(V_{N,k}^{r,2}\right) \Big) \Big] dr .
	\end{align*}
	
	\noindent Since $P_{1,2}$ is a matrix of rank $N$, we now fix $D = \max(1,\sup_{u,v } \norm{L_j^{u,v}})$, we have 
	$$ \left| \cov_k\left( \tr_{kN}(A_s C_sP_{1,2}), \tr_{kN}(B_s) \right) \right| \leq \frac{s}{k} \deg(AC) \deg(B)\ D^{\deg(ABC)} .$$
	
	\noindent We now assume that $Q$ has degree at most $n$, then $ \left| \cov_k\left( \tr_{kN}(A_s C_sP_{1,2}), \tr_{kN}(B_s) \right) \right| \leq \frac{s}{k} n^2 D^{n}$. Thus we have,
	\begin{align*}
	\left| \E_k\left[\tr_{kN}\left( Q_t P_{1,2}\right)\right] \right| \leq \frac{n^4 t^2 D^n}{k^2N} \\
	+ \frac{1}{kN} \int_{0}^t e^{- \frac{|Q|_B}{2} (t-s)} \Bigg( & \sum_{i\leq p,\  Q = A U_i B U_i C} \left| \E_k\left[\tr_{kN}\left(A_s U^{kN}_{i,s} C_sP_{1,2}\right)\right] \E_k\left[\tr_{kN}\left(B_s U^{kN}_{i,s}\right)\right] \right| \\
	& \sum_{i\leq p,\  Q = A U_i^* B U_i^* C} \left| \E_k\left[\tr_{kN}\left(A_s {U^{kN}_{i,s}}^* C_tP_{1,2}\right)\right] \E_k\left[\tr_{kN}\left(B_s {U^{kN}_{i,s}}^*\right)\right] \right| \\
	&\sum_{i\leq p,\ Q=A U_i B U_i^* C } \left|\E_k\left[\tr_{kN}(A_s C_sP_{1,2})\right] \E_k\left[\tr_{kN}(B_s)\right] \right| \\
	&\sum_{i\leq p,\ Q=A U_i^* B U_i C } \left|\E_k\left[\tr_{kN}(A_s C_sP_{1,2})\right] \E_k\left[\tr_{kN}(B_s)\right] \right| \Bigg) ds .
	\end{align*}
	
	\noindent This means that,	
	\begin{align*}
	\left| \E_k\left[\tr_{kN}\left( Q_t P_{1,2}\right)\right] \right| &\leq \frac{n^4 t^2 D^n}{k^2N}	+ \int_{0}^t e^{- \frac{|Q|_B}{2} (t-s)} \Bigg( \sum_{1\leq i\leq p,\  Q = A U_i B U_i C} c^t_{\deg(AC)+1} D^{\deg(B)+1} \\
	&\quad\quad\quad\quad\quad\quad\quad\quad\quad\quad\quad\quad\quad\ + \sum_{1\leq i\leq p,\  Q = A U_i^* B U_i^* C} c^t_{\deg(AC)+1} D^{\deg(B)+1} \\
	&\quad\quad\quad\quad\quad\quad\quad\quad\quad\quad\quad\quad\quad\ + \sum_{1\leq i\leq p,\ Q=A U_i B U_i^* C } c^t_{\deg(AC)} D^{\deg(B)} \\
	&\quad\quad\quad\quad\quad\quad\quad\quad\quad\quad\quad\quad\quad\ + \sum_{1\leq i\leq p,\ Q=A U_i^* B U_i C } c^t_{\deg(AC)} D^{\deg(B)} \Bigg) ds \\
	&\leq \frac{n^4 t^2 D^n}{k^2N}	+ \int_{0}^t |Q|_B e^{- \frac{|Q|_B}{2} s} ds \sum_{0\leq d\leq n-1} D^d c_{n-1-d}^t \\
	&\leq \frac{n^4 t^2 D^n}{k^2N}	+ 2 \sum_{0\leq d\leq n-1} D^d c_{n-1-d}^t .
	\end{align*}
	
	\noindent Hence, for any $n\geq 1$,
	$$ c_n^t \leq \frac{n^4 t^2 D^n}{k^2N}	+ 2 \sum_{0\leq d\leq n-1} D^d c_{n-1-d}^t .$$
	
	\noindent Since we are taking the trace of $L(V_{N,k}^s) P_{1,2}$ with $P_{1,2} = I_N \otimes E_{1,2}$, we have $c_0 = 0$. We fix $s:a\mapsto \sum_{n\geq 0} \frac{n^4 t^2 (aD)^n}{N}$, thus for $a$ small enough, 
	\begin{align*}
	f_k^t(a) &\leq \frac{s(a)}{k^2} + 2 \sum_{n\geq 1} \left(\sum_{0\leq d\leq n-1} D^d c_{n-1-d}^t \right) a^n \\
	&\leq \frac{s(a)}{k^2} + \frac{2a}{1-aD} f_k^t(a)
	\end{align*}
	
	\noindent Thus for $a$ small enough, $f_k^t(a) \leq 2 s(a) k^{-2}$. Which means that $f_k^t(a) = \mathcal{O}(k^{-2})$, hence the conclusion.
	
\end{proof}

We can now prove the following intermediary lemma that will allow us to derive Lemma \ref{2optilll}. This lemma is the only one where the law of $U^N$ actually plays an important part. To be more precise, we use the invariance of the law of a Haar unitary matrix by multiplication by a deterministic unitary matrix.

\begin{lemma}
	\label{2grcalc}
	We define $Y_t^k$ as in Proposition \ref{2limite_trace}, we set
	\begin{itemize}
		\item $P_{l,l'} = I_N\otimes E_{l,l'}\otimes I_M$,
		\item $q = A\ e^{\i y P}\ B$.
	\end{itemize}
	
	\noindent Then for every $M,N\in\N$, $t\in\R^+$ and $y\in\R$,
	
	\begin{align*}
	S_{t,y}^N(A,B) = \lim_{k\to\infty} - \frac{1}{k N^2} \E\Bigg[\tau_M\Bigg( \sum_{1\leq l,l'\leq k} &\E_k\Big[ (\tr_{kN}\otimes I_{M})^{\bigotimes 2}\left(\delta_i q(Y_t^k) \times P_{l',l}\otimes P_{l,l'} \right) \Big] \\
	& - \E_k[\tr_{kN}\otimes I_{M}]^{\bigotimes 2}\left(\delta_i q(Y_t^k) \times P_{l',l}\otimes P_{l,l'} \right) \Bigg)\Bigg]
	\end{align*}
	where thanks to Proposition \ref{2duhamel}, we set
	\begin{equation}
	\label{2seriveqtr}
	\delta_i q = \delta_i A\  e^{\i y P}\ B + \i y A\ \int_0^1  e^{\i \alpha y P}\ \delta_i P\ e^{\i (1-\alpha) y P}\ B d\alpha + A\  e^{\i y P}\ \delta_i B .
	\end{equation}
	
\end{lemma}

\begin{proof}
	
	Since all of our random variables are unitary matrices, thanks to Proposition \ref{2limite_trace} and the dominated convergence theorem,
	\begin{equation}
	\label{2premiertr}
	S_{t,y}^N(A,B) = \lim_{k\to\infty} \E\left[ \tau_M\left( \E_k[\tau_{kN}\otimes I_M] \bigotimes \E_k[\tau_{kN}\otimes I_M]\left(\delta_i\left( A\ e^{\i y P}\ B \right) \left(Y_t^k\right)\right)\right) \right] ,
	\end{equation}

	\noindent where $\E_k[\tau_{kN}\otimes I_M] \bigotimes \E_k[\tau_{kN}\otimes I_M]\left( A\otimes B \left(Y_t^k\right)\right) = \E_k[\tau_{kN}\otimes I_M (A(Y_t^k))]\E_k[\tau_{kN}\otimes I_M (B(Y_t^k))]$. Since given $V\in \U_N$,  $(U^{kN}_{t,1},U^N_1\otimes I_k,\dots,U^{kN}_{t,p},U^N_p\otimes I_k)$ has the same law as $((V^*\otimes I_k) U^{kN}_{t,1} (V\otimes I_k),(U^N_1 V)\otimes I_k,U^{kN}_{t,2},U^N_2\otimes I_k,\dots,U^{kN}_{t,p},U^N_p\otimes I_k)$, we have
	
	\begin{align*}
	\E[ q(Y_t^k) ] = \E\Big[ \widetilde{q}\Big( &(U^N_1\otimes I_k\ U^{kN}_{t,1})\otimes I_M\ (V\otimes I_{kM}), \\
	&(U^N_2\otimes I_k\ U_{t,2}^{kN})\otimes I_M,\dots,(U^N_p\otimes I_k\ U_{t,p}^{kN})\otimes I_M, Z^{NM}\otimes I_k\Big) \Big] .
	\end{align*}
	
	\noindent Hence let $H$ be an skew-Hermitian matrix, then for any $s\in \R$, $e^{sH}\in \U_N$, thus by taking $V$ this matrix and differentiating with respect to $s$ we get that, $\E\Big[ \delta_1 q(Y_t^k) \# (H\otimes I_{kM}) \Big] =0$. And similarly we get that for any $i$,
	$$\E\Big[ \delta_i q(Y_t^k) \# (H\otimes I_{kM}) \Big] =0 .$$
	
	\noindent Since every matrix is a linear combination of skew-Hermitian matrices (indeed if $A\in\M_N(\C)$, then $2A = (A-A^*) + \i \times (-\i)(A^*+A)\ $), this is true for any matrix $H\in \M_N(\C)$, and thus for any $i$,
	\begin{equation}
	\label{2SDeq}
	\E\Big[ (\tr_N\otimes I_{kM})^{\bigotimes 2}\left(\delta_i q(Y_t^k)\right) \Big] = \sum_{ 1\leq r,s \leq N} g_r^*\otimes I_{kM}\ \E\Big[ \delta_i q(Y_t^k) \# (E_{r,s}\otimes I_{kM}) \Big] g_s\otimes I_{kM} = 0
	\end{equation}
	
	\noindent Let $S,T\in \M_{NkM}(\C)$ be deterministic matrices, then
	\begin{align*}
	& \tr_{k}\otimes I_M\left( (\tr_N\otimes I_{kM})^{\bigotimes 2}\left( S\otimes T \right) \right) \\
	&= \sum_{1\leq m,n\leq N} \tr_{Nk}\otimes I_M\left(S\ E_{m,n}\otimes I_{kM}\ T\ E_{n,m}\otimes I_{kM}\right)  \\
	&= \sum_{1\leq l,l'\leq k} \sum_{1\leq m\leq N} g_m^*\otimes f_l^*\otimes I_M\ S\ g_m\otimes f_{l'}\otimes I_M \sum_{1\leq n\leq N} g_n^*\otimes f_{l'}^*\otimes I_M\ T\ g_n\otimes f_l\otimes I_M  \\
	&= \sum_{1\leq l,l'\leq k} \tr_N\otimes I_M( I_N\otimes f_l^*\otimes I_M\ S\ I_N\otimes f_{l'}\otimes I_M) \tr_N\otimes I_M(I_N\otimes f_{l'}^*\otimes I_M\ T\ I_N\otimes f_l\otimes I_M) \\
	&= \sum_{1\leq l,l'\leq k} \tr_{kN}\otimes I_M\big( S\ I_N\otimes E_{l',l}\otimes I_M\big) \tr_{kN}\otimes I_M\big( T\ I_N\otimes E_{l,l'}\otimes I_M\big) .
	\end{align*}
	
	\noindent Thus by using equation \eqref{2SDeq}, we have for any $i$,
	\begin{equation*}
	\sum_{1\leq l,l'\leq k} \E\Big[ (\tr_{kN}\otimes I_{M})^{\bigotimes 2}\left(\delta_i q(Y_t^k) \times P_{l',l}\otimes P_{l,l'} \right) \Big] = 0 .
	\end{equation*}
	
	\noindent And consequently,
	\begin{align}
	\label{2deuxiemetr}
	&\sum_{1\leq l,l'\leq k} \E\Big[ (\tr_{kN}\otimes I_{M})^{\bigotimes 2}\left(\delta_i q(Y_t^k) \times P_{l',l}\otimes P_{l,l'} \right) \Big] \nonumber \\
	&\quad\quad - \E\left[ \E_k[\tr_{kN}\otimes I_{M}]^{\bigotimes 2}\left(\delta_i q(Y_t^k) \times P_{l',l}\otimes P_{l,l'} \right) \right] \\
	&= - \sum_{1\leq l,l'\leq k} \E\left[ \E_k[\tr_{kN}\otimes I_{M}]^{\bigotimes 2}\left(\delta_i q(Y_t^k) \times P_{l',l}\otimes P_{l,l'} \right) \right] . \nonumber
	\end{align}
	
	\noindent Let $V,W\in \M_k(\C)$ be permutation matrices. Since $I_{NM}\otimes V$ commutes with $Z^{NM}\otimes I_k$ and $U^N\otimes I_{kM}$, and that the law of $U^{kN}_t$ is invariant by conjugation by a unitary matrix, it follows that the law of every matrix of $Y_t^k$ is invariant by conjugation by $I_{NM}\otimes V$ or $I_{NM}\otimes W$. Thus,
	$$ \E_k[\tr_{kN}\otimes I_{M}]^{\bigotimes 2}\left(\delta_i q(Y_t^k) \times P_{l',l}\otimes P_{l,l'} \right) = \E_k[\tr_{kN}\otimes I_{M}]^{\bigotimes 2}\left(\delta_i q(Y_t^k) \times VP_{l',l}V^*\otimes WP_{l,l'}W^* \right) . $$
	
	\noindent Thus by using well-chosen matrices, we get
	\begin{itemize}
		\item if $l= l'$, $$\E_k[\tr_{kN}\otimes I_{M}]^{\bigotimes 2}\left(\delta_i q(Y_t^k) \times P_{l',l}\otimes P_{l,l'} \right) = \E_k[\tr_{kN}\otimes I_{M}]^{\bigotimes 2}\left(\delta_i q(Y_t^k) \times P_{1,1}\otimes P_{1,1} \right),$$
		\item if $l\neq l'$, $$\E_k[\tr_{kN}\otimes I_{M}]^{\bigotimes 2}\left(\delta_i q(Y_t^k) \times P_{l',l}\otimes P_{l,l'} \right) = \E_k[\tr_{kN}\otimes I_{M}]^{\bigotimes 2}\left(\delta_i q(Y_t^k) \times P_{1,2}\otimes P_{1,2} \right).$$
	\end{itemize}
	
	\noindent Consequently, we have that
	\begin{itemize}
		\item equation \eqref{2deuxiemetr} is equal to
		\begin{align*}
		&\sum_{1\leq l,l'\leq k} \E\Big[ (\tr_{kN}\otimes I_{M})^{\bigotimes 2}\left(\delta_i q(Y_t^k) \times P_{l',l}\otimes P_{l,l'} \right) \Big] \\
		&\quad\quad - \E\left[ \E_k[\tr_{kN}\otimes I_{M}]^{\bigotimes 2}\left(\delta_i q(Y_t^k) \times P_{l',l}\otimes P_{l,l'} \right) \right] \\
		&= - k \E\left[\E_k[\tr_{kN}\otimes I_{M}]^{\bigotimes 2}\left(\delta_i q(Y_t^k) \times P_{1,1}\otimes P_{1,1} \right)\right] \\
		&\quad - k(k-1) \E\left[\E_k[\tr_{kN}\otimes I_{M}]^{\bigotimes 2}\left(\delta_i q(Y_t^k) \times P_{1,2}\otimes P_{1,2} \right) \right] .
		\end{align*}
		\item Whereas the quantity inside the trace $\tau_M$ in equation \eqref{2premiertr} is equal to
		\begin{align*}
		&\E\left[\E_k[\tau_{kN}\otimes I_M] \bigotimes \E_k[\tau_{kN}\otimes I_M]\left(\delta_iq \left(Y_t^k\right)\right)\right] \\
		&= \frac{1}{(kN)^2} \sum_{ 1\leq l,l' \leq k} \E\left[\E_k[\tr_{kN}\otimes I_M]^{\bigotimes 2}\left(\delta_iq \left(Y_t^k\right) \times P_{l,l}\otimes P_{l',l'}\right)\right] \\
		&= \frac{1}{N^2} \E\left[\E_k[\tr_{kN}\otimes I_M]^{\bigotimes 2}\left(\delta_iq \left(Y_t^k\right) \times P_{1,1}\otimes P_{1,1}\right)\right] .
		\end{align*}
	\end{itemize}
	
	\noindent Thus, we have
	\begin{align*}
	S_{t,y}^N(A,B) = \lim_{k\to\infty} - &\frac{1}{k N^2} \tau_M\Bigg( \sum_{1\leq l,l'\leq k} \E\Big[ (\tr_{kN}\otimes I_{M})^{\bigotimes 2}\left(\delta_i q(Y_t^k) \times P_{l',l}\otimes P_{l,l'} \right) \Big] \\
	&\quad\quad\quad\quad\quad\quad\quad\quad - \E\left[\E_k[\tr_{kN}\otimes I_{M}]^{\bigotimes 2}\left(\delta_i q(Y_t^k) \times P_{l',l}\otimes P_{l,l'} \right)\right] \Bigg) \\
	-& \frac{k-1}{N^2} \E\left[ \tau_M\Bigg( \E_k[\tr_{kN}\otimes I_{M}]^{\bigotimes 2}\left(\delta_i q(Y_t^k) \times P_{1,2}\otimes P_{1,2} \right) \Bigg) \right] .
	\end{align*}
	
	\noindent Thanks to Lemma \ref{2nondiag} and Proposition \ref{2duhamel}, the last term converges towards $0$, which gives the expected formula.
\end{proof}

Lemma \ref{2grcalc} makes a covariance appears. Thus it is natural to want to use Corollary \ref{2concentrationmulitdim} to get an upper bound of $S_{t,y}^N(A,B)$, explicit in all of its parameters.

\begin{lemma}
	\label{2optilll}
	There exists a polynomial $L_P\in\R^+[X]$ such that for any $t,y,N,M,Z^{NM}$, 
	\begin{equation}
	\label{2optiN}
	|S_{t,y}^N(A,B)|  \leq L_P\left( \norm{Z^{NM}} \right) \frac{M^2}{N^2} (1+|y|^3)\ t.
	\end{equation}
	Besides if $Z^{NM} = (I_N\otimes Y_1^M,\dots ,I_N\otimes Y_q^M)$ and that these matrices commute, then we have the same inequality without the $M^2$.
	
\end{lemma}

\begin{proof}
	As mentioned in equation \eqref{2seriveqtr}, we have
	\begin{equation*}
	\delta_i q = \delta_i A\  e^{\i y P}\ B + \i y A\ \int_0^1  e^{\i \alpha y P}\ \delta_i P\ e^{\i (1-\alpha) y P}\ B d\alpha + A\  e^{\i y P}\ \delta_i B .
	\end{equation*}
	
	\noindent Consequently, we set $q_1 = A_1 e^{\i \alpha y P} B_1$ and $q_2 = A_2 e^{\i (1-\alpha) y P} B_2$ where $A_1,B_1,A_2,B_2\in\PP_{d}$ are monomials, then thanks to equation \eqref{2extension} and Proposition \ref{2duhamel}, we can use Corollary \ref{2concentrationmulitdim} even though $q_1$ and $q_2$ are not exactly polynomials and we get that
	\begin{align*}
	&\tau_M\Big(\E_k\Big[ (\tr_{kN}\otimes I_{M})^{\bigotimes 2}\left(q_1(Y_t^k)P_{l',l}\otimes q_2(Y_t^k)P_{l,l'} \right) \Big] \\
	&\quad\quad\quad\quad\quad\quad\quad\quad\quad\quad\quad\quad- \E_k[\tr_{kN}\otimes I_{M}]^{\bigotimes 2}\left(q_1(Y_t^k)P_{l',l}\otimes q_2(Y_t^k)P_{l,l'} \right) \Big) \\
	&= \sum_{ 1\leq j \leq p}\int_{0}^{t} \E_k\Big[ \tau_{kN}\otimes\tau_M\Big( h\circ\delta_j(\widetilde{q_1} P_{l',l})(U^N\otimes I_k V_{t-s}^{kN}U_s^{kN},Z^{NM}\otimes I_M) \\
	&\quad\quad\quad\quad\quad\quad\quad\quad\quad\quad\quad \times h\circ\delta_j(\widetilde{q_2} P_{l',l})(U^N\otimes I_k W_{t-s}^{kN}U_s^{kN},Z^{NM}\otimes I_M)\Big) \Big]\ ds,
	\end{align*}
	where $(V_{s}^{kN})_{s\in\R}$ and $(W_{s}^{kN})_{s\in\R}$ are $p$-tuples of independent unitary Brownian motions of size $kN$. Thus thanks to Lemma \ref{2normineq}, we get that there exist a polynomial $L_{A_1,B_1,A_2,B_2,P}$ such that 
	\begin{align*}
	&\Big| \tau_M\Big( \E_k\Big[ (\tr_{kN}\otimes I_{M})^{\bigotimes 2}\left(q_1(Y_t^k)P_{l',l}\otimes q_2(Y_t^k)P_{l,l'} \right) \Big] \\
	&\quad\quad\quad\quad\quad\quad\quad\quad\quad\quad\quad\quad  - \E_k[\tr_{kN}\otimes I_{M}]^{\bigotimes 2}\left(q_1(Y_t^k)P_{l',l}\otimes q_2(Y_t^k)P_{l,l'} \right)\Big) \Big| \\
	&\leq L_{A_1,B_1,A_2,B_2,P}\left(\norm{Z^{NM}}\right) \times (1+y^2) \frac{M^2t}{k},
	\end{align*}
	where we used the fact that $P_{l',l}$ has rank $NM$ and that the renormalized trace of a matrix of rank $NM$ in $\M_{kN}(\C)\otimes\M_M(\C)$ is smaller than its norm renormalized by $k$.
	
	Since this upper bound does not depend on $\alpha$, it remains true if we integrate with respect to $\alpha$ from $0$ to $1$. But then $\delta_i q$ is a finite linear combination of such terms. Consequently, one gets that there exists a polynomial $L_P$ such that for any $k$, 
	\begin{align*}
	&\Bigg| \frac{1}{k N^2} \tau_M\Bigg( \sum_{1\leq l,l'\leq k} \E\Big[ (\tr_{kN}\otimes I_{M})^{\bigotimes 2}\left(\delta_i q(Y_t^k) \times P_{l',l}\otimes P_{l,l'} \right) \Big] \\
	&\quad\quad\quad\quad\quad\quad\quad\quad - \E\left[\E_k[\tr_{kN}\otimes I_{M}]^{\bigotimes 2}\left(\delta_i q(Y_t^k) \times P_{l',l}\otimes P_{l,l'} \right)\right] \Bigg) \Bigg| \\
	&\leq L_{P}\left(\norm{Z^{NM}}\right) \times (1+|y|^3) \frac{M^2t}{N^2}.
	\end{align*}
	Finally, when the matrices $Z^{NM}$ commute, as specified in Lemma \ref{2concentrationmulitdim}, we have the same proof where we replaced $h\circ \delta_j$ by $\D_j$, hence we do not need to use Lemma \ref{2normineq} and hence we have the same inequality without the $M^2$. Finally we get the conclusion thanks to Lemma \ref{2grcalc}.
\end{proof}

We now have the tools to prove Theorem \ref{2imp2}.

\begin{proof}[Proof of Theorem \ref{2imp2}]
	
	Thanks to Lemma \ref{2fourier}, and since $\D_i\ e^{\i y P} = \i y\ \delta_i P\ \widetilde{\#}\ e^{\i y P}$, there exist a family of monomials $(A_k,B_k)_k$ and a constant $C$ which only depends on $P$ such that,
	\begin{align*}
	\Bigg| &\E\left[\frac{1}{MN}\tr_{MN}\Big(f\left(\widetilde{P}\left(U^N\otimes I_M,Z^{NM}\right)\right)\Big) - \tau_N\otimes\tau_M\Big(f\left(\widetilde{P}\left(U^Nu_T\otimes I_M,Z^{NM}\right)\right)\Big) \right] \Bigg| \\
	&\leq C \sum_{k} \int |y| \int_0^T \left| S_{t,y}^N(A_k,B_k) \right| dt\ d|\mu|(y) .
	\end{align*}
	
	\noindent Thanks to equation \eqref{2optiN}, we get that for some polynomial $L_P$,
	\begin{align*}
	\Bigg| &\E\left[\frac{1}{MN}\tr_{MN}\Big(f\left(\widetilde{P}\left(U^N\otimes I_M,Z^{NM}\right)\right)\Big) - \tau_N\otimes\tau_M\Big(f\left(\widetilde{P}\left(U^Nu_T\otimes I_M,Z^{NM}\right)\right)\Big) \right] \Bigg| \\
	&\leq T^2\  \frac{M^2}{N^2} L_P\left(\norm{Z^{NM}}\right)\times \int_{\R} (|y|+y^4)\  d|\mu|(y) \ .
	\end{align*}
	
	\noindent And besides if $Z^{NM} = (I_N\otimes Y_1^M,\dots ,I_N\otimes Y_q^M)$ and that these matrices commute, we have the same inequality without the $M^2$. Finally, thanks to Proposition \ref{2Browniantounitary}, thanks to Duhamel's formula \eqref{2Duhmeqned} we can find a polynomial $L_P'$ such that
	\begin{align*}
	&\Bigg| \tau_N\otimes\tau_M\Big(e^{\i y\widetilde{P}\left(u\otimes I_M,Z^{NM}\right)}\Big) - \tau_N\otimes\tau_M\Big(e^{\i y\widetilde{P}\left(U^Nu_T\otimes I_M,Z^{NM}\right)}\Big) \Bigg| \\
	&= \Bigg| \tau_N\otimes\tau_M\Big(e^{\i y\widetilde{P}\left(U^N u\otimes I_M,Z^{NM}\right)}\Big) - \tau_N\otimes\tau_M\Big(e^{\i y\widetilde{P}\left(U^Nu_T\otimes I_M,Z^{NM}\right)}\Big) \Bigg| \\
	&\leq e^{-T/2} L_P'\left(\norm{Z^{NM}}\right)\times |y| .
	\end{align*}
	
	\noindent Hence the conclusion by fixing $T = 4\ln(N)$.
	
\end{proof}

We can finally prove Theorem \ref{2imp1}.

\begin{proof}[Proof of Theorem \ref{2imp1}]
	
	We want to use Theorem \ref{2imp2}. To do so we would like to take the Fourier transform of $f$ and use Fourier inversion formula. However we did not assume that $f$ is integrable. Thus the first step of the proof is to show that we can assume that $f$ has compact support. Since $U^N$ and $u$ are unitaries, there exists a polynomial $H\in\R^+[X]$ which only depends on $P$ such that $ \norm{\widetilde{P}\left(U^N\otimes I_M,Z^{NM}\right)} \leq H\left( \norm{Z^{NM}} \right)$. Consequently since we also have that $\norm{\widetilde{P}(u\otimes I_M,Z^{NM})} \leq H\left( \norm{Z^{NM}} \right)$, we can replace $f$ by $fg$ where $g$ is a $\mathcal{C}^{\infty}$ function which takes value in $[0,1]$, takes value $1$ in $[-H\left( \norm{Z^{NM}} \right),H\left( \norm{Z^{NM}} \right)]$ and $0$ outside of $[-H\left( \norm{Z^{NM}} \right) -1,H\left( \norm{Z^{NM}} \right) +1]$. Since $f$ can be differentiated six times, we can take its Fourier transform and then invert it so that with the convention $ \hat{f}(y) = \frac{1}{2\pi} \int_{\R} f(x) e^{-\i xy} dx$, we have
	
	$$ \forall x\in\R,\quad f(x) = \int_{\R} e^{\i xy} \widehat{f}(y)\ dy . $$
	
	\noindent Besides since if $f$ has compact support bounded by $H\left( \norm{Z^{NM}} \right) +1$, then 
	$$\norm{\hat{f}}_{\infty} \leq 2\left(H\left( \norm{Z^{NM}} \right) +1\right) \norm{f}_{\infty} ,$$
	we get that
	
	\begin{align*}
	\int_{\R} (|y|+y^4)\ \left| \widehat{f}(y) \right|\ dy &\leq \int_{\R} \frac{|y|+|y|^3+y^4+y^6}{1+y^2}\ \left| \widehat{f}(y) \right|\ dy \\
	&\leq \bigintss_{\R} \frac{\left| \widehat{(f)^{(1)}}(y) \right| + \left| \widehat{(f)^{(3)}}(y) \right| + \left| \widehat{(f)^{(4)}}(y) \right| + \left| \widehat{(f)^{(6)}}(y) \right|}{1+y^2}\ dy \\
	&\leq 2 \left( H\left(\norm{Z^{NM}}\right) + 1\right) \norm{f}_{\mathcal{C}^6} \int_{\R} \frac{1}{1+y^2}\ dy \\
	&\leq 2\pi \left( H\left(\norm{Z^{NM}}\right) + 1\right) \norm{f}_{\mathcal{C}^6} ,
	\end{align*}
	
	\noindent Hence it satisfies the hypothesis of Theorem \ref{2imp2} with $\mu(dy) = \widehat{f}(y) dy$, thus we have
	
	\begin{align*}
	&\left| \E\left[\frac{1}{MN}\tr\Big(f\left(\widetilde{P}\left( U^N\otimes I_M,Z^{NM}\right)\right)\Big)\right] - \tau\Big(f\left(\widetilde{P}\left( u\otimes I_M,Z^{NM}\right)\right)\Big) \right| \nonumber \\
	&\leq \frac{M^2\ln^2(N)}{N^2} L_P\left(\norm{Z^{NM}}\right) \int_{\R} (|y|+y^4)\ \left| \widehat{f}(y) \right|\ dy \\
	&\leq \frac{M^2\ln^2(N)}{N^2}\ \times 2\pi L_P\left(\norm{Z^{NM}}\right) \left( H\left(\norm{Z^{NM}}\right) + 1\right) \norm{f}_{\mathcal{C}^6} . \\
	\end{align*}
	
	\noindent And finally, if $Z^{NM} = (I_N\otimes Y_1^M,\dots ,I_N\otimes Y_q^M)$ and that these matrices commute, then we have the same inequality without the $M^2$ as specified in Theorem \ref{2imp2}. 
	
\end{proof}

\vspace*{1cm}

\section{Proof of Corollaries}
\label{2proofcoro}

\subsection{Proof of Corollary \ref{2stieljes}}

We could directly apply Theorem \ref{2imp1} to $f_z : x\to (z-x)^{-1}$, however for $z$ such that $\Im z$ is small, we have $\norm{f}_{\mathbb{C}^6} = O\left((\Im z)^{-7} \right) $ when we want $O\left((\Im z)^{-5} \right) $ instead. Since $P$ is self-adjoint, $\overline{G_{P}(z)} = G_{P}(\overline{z})$, thus we can assume that $\Im z < 0$, but then

$$ f_z(x) = \int_{0}^{\infty} e^{\i x y}\ (\i e^{-\i y z})\ dy .$$

\noindent Consequently with $\mu_z(dy) = \i e^{-\i y z}\ dy $, we have

\begin{align*}
\int_{0}^{\infty} (y+y^4)\ d|\mu_z|(y) &= \frac{1}{|\Im z|^2} + \frac{24}{|\Im_z|^5} .
\end{align*}

\noindent Thus by applying Theorem \ref{2imp2} with $Z^{NM} = \left(I_N\otimes Y_1^M,\dots,I_N\otimes Y_p^M\right)$, $P$ and $f_z$, we have

$$ \left| \E\left[ G_{P(U^N\otimes I_M, I_N\otimes Y^M)}(z) \right] - G_{P(u\otimes I_M, 1\otimes Y^M)}(z) \right| \leq \frac{M^2 \ln^2(N)}{N^2} L_P\left(\norm{Z^{NM}}\right) \int_{\R} (1+y^4)\  d|\mu_z|(y) .$$

\noindent Now since $\norm{Z^{NM}} = \norm{Y^M}$ which does not depend on $N$, we finally have
\begin{align*}
&\left| \E\left[ G_{P(U^N\otimes I_M, I_N\otimes Y^M)}(z) \right] - G_{P(u\otimes I_M, 1\otimes Y^M)}(z) \right| \leq \frac{M^2\ln^2(N)}{N^2} L_P\left(\norm{Y^M}\right)  \left(\frac{1}{|\Im z|^2} + \frac{24}{|\Im_z|^5}\right) .
\end{align*}

\subsection{Proof of Corollary \ref{2LU}}

Let $f:\R\to\R$ be a Lipschitz function uniformly bounded by $1$ and with Lipschitz constant at most $1$, we want to find an upper bound on

\begin{equation}
\label{2initi}
\Bigg| \E\left[\frac{1}{MN}\tr_{NM}\Big(f\left(P\left(U^N\otimes I_M,I_N\otimes Y_M\right)\right)\Big)\right] - \tau\otimes\tau_M\Big(f\left(P\left(u\otimes I_M,1\otimes Y_M\right)\right)\Big) \Bigg| .
\end{equation}

\noindent Firstly, since $U^N$ are unitary matrices, we can assume that the support of $f$ is bounded by a constant $S = H(\norm{Y^M})$ for some polynomial $H\in\R^+[X]$ independent of everything. However we cannot apply directly Theorem \ref{2imp1} since $f$ is not regular enough. In order to deal with this issue we use the convolution with gaussian random variable, thus let $G$ be a centered gaussian random variable, we set

$$ f_{\varepsilon} : x\to \E[f(x+\varepsilon G)] .$$

\noindent Since $f$ has Lipschitz constant $1$, we have for any $x\in\R$,

$$ \left| \E[f(x+\varepsilon G)] -f(x)\right| \leq \varepsilon .$$

\noindent Since $f_{\varepsilon}$ is regular enough we could now apply Theorem \ref{2imp1}, however we a get better result by using Theorem \ref{2imp2}. Indeed we have

\begin{align*}
f_{\varepsilon}(x) &= \frac{1}{\sqrt{2\pi}} \int_{\R} f(x+\varepsilon y) e^{-y^2/2}\  dy \\
&= \frac{1}{\sqrt{2\pi}} \int_{\R} f(y) \frac{e^{-\frac{(x-y)^2}{2\varepsilon^2}}}{\varepsilon}\  dy \\
&= \frac{1}{2\pi} \int_{\R} f(y) \int_{\R} e^{\i (x-y) u} e^{-(u\varepsilon)^2/2}\ du\ dy .
\end{align*}

\noindent Since the support of $f$ is bounded, we can apply Fubini's theorem:

\begin{align*}
f_{\varepsilon}(x) &= \frac{1}{2\pi} \int_{\R} e^{\i ux} \int_{\R} f(y) e^{- \i y u}  dy\  e^{-(u\varepsilon)^2/2}\  du .
\end{align*}

\noindent And so with the convention $ \hat{h}(u) = \frac{1}{2\pi} \int_{\R} h(y) e^{-\i uy} dy$, we have

$$ f_{\varepsilon}(x) = \int_{\R} e^{\i ux} \hat{f}(u)  e^{-(u\varepsilon)^2/2} du .$$

\noindent Thus if we set $\mu_{\varepsilon}(dy) = \hat{f}(y)  e^{-(y\varepsilon)^2/2} dy$, then, since $\norm{f}_{\infty} \leq 1$,

$$ \int_{\R} (1+y^4) d|\mu_{\varepsilon}|(y) \leq 2S \int_{\R} (1+y^4)e^{-y^2/2}\ dy\ \varepsilon^{-5} . $$

\noindent Consequently we can apply Theorem \ref{2imp2} with $f_{\varepsilon}$ and since $\norm{f-f_{\varepsilon}}_{\infty}\leq \varepsilon$, there exists a polynomial $R_P$ such that \eqref{2initi} can be bounded by

$$ 2\varepsilon + R_P\left(\norm{Y^M}\right) \frac{M^2 \ln^2(N)}{N^2 \varepsilon^5} .$$

\noindent Thus we can now fix $\varepsilon = (N^{-1} \ln(N))^{1/3}$ and we get that for any $f$ Lipschitz function uniformly bounded by $1$ and with Lipschitz constant at most $1$, \eqref{2initi} can be bounded by

$$ 2 R_P\left(\norm{Y^M}\right) M^2 \left( \frac{\ln N}{N} \right)^{1/3} . $$

\subsection{Proof of Theorem \ref{2strongconv}}

Firstly, we need to set the operator norm of tensor of $\CC^*$-algebras we will work with. When writing the proof it appears that it is the minimal tensor product as defined in \ref{2mini}. The following two lemmas were used in \cite{un}, see Lemma 4.1.8 from \cite{ozabr} for a proof of the first one and Lemma 4.3 from \cite{un} for the second one. In order to learn more about the second lemma, especially how to weaken the hypothesis, we refer to \cite{pisier}.

\begin{lemma}
	\label{2faith}
	Let $(\A,\tau_{\A})$ and $(\B,\tau_{\B})$ be $\CC^*$-algebras with faithful traces, then $\tau_{\A}\otimes\tau_{\B}$ extends uniquely to a faithful trace $\tau_{\A}\otimes_{\min}\tau_{\B}$ on $\A\otimes_{\min}\B$. 
\end{lemma}

\begin{lemma}
	\label{2tensorconv}
	Let $\mathcal{C}$ be an exact $\mathcal{C}^*$-algebra endowed with a faithful state $\tau_{\CC}$, let $Y^N \in \mathcal{A}_N$ be a sequence of family of noncommutative random variable in a $\mathcal{C}^*$-algebra $\mathcal{A}_N$ which converges strongly towards a family $Y$ in a $\mathcal{C}^*$-algebra $\mathcal{A}$ endowed with a faithful state $\tau_{\A}$. Let $S\in \mathcal{C}$ be a family of noncommutative random variable, then the family $(S\otimes 1, 1\otimes Y^N)$ converges strongly in distribution towards the family $(S\otimes 1, 1\otimes Y)$.
\end{lemma}

In order to prove Theorem \ref{2strongconv} we use well-known concentration properties of unitary Haar matrices coupled with 
an estimation of the expectation, let us begin by stating the concentration properties that we will use.

\begin{prop}
	\label{2concentr}
	Let $f$ be a continuous function on $\U_N^p$, such that for any $X,Y\in \U_N^p$,
	$$ |f(X)-f(Y)| \leq C \sum_i \tr_N\left((X_i-Y_i)(X_i-Y_i)^*\right)^{1/2} .$$
	Then if $W$ is a family of $p$ independent random matrices distributed according to the Haar measure on $\mathbb{SU}_N$, and $U$ a family of $p$ independent unitary Haar matrices of size $N$ independent from $W$, we have,
	$$ \P\left( \left| f(U) - \E_W\left[f(W U)\right]  \right| \geq \delta \right) \leq 4p\ e^{- \left(\frac{\delta}{2 p C}\right)^2 N} .$$
\end{prop}

\begin{proof}
	We want to use Corollary 4.4.28 from \cite{alice}, in order to do so let us first assume that $f$ takes real values. We then set,
	$$ f^i_{U_{i+1},\dots,U_p} : U_i \to \E_{W_1,\dots,W_{i-1}}\left[f(W_1 U_1,\dots, W_{i-1} U_{i-1}, U_i, U_{i+1},\dots,U_p )\right] .$$
	
	\noindent Thus,
	$$ f(U) - \E_W\left[f(W U)\right] = \sum\limits_{ 1\leq i \leq p} f^i_{U_{i+1},\dots,U_p}(U_i) - \E_{W_i}\left[ f^i_{U_{i+1},\dots,U_p}(W_i U_i) \right] .$$
	
	\noindent Besides for any $U_i,V_i$, we have that 
	$$ | f^i_{U_{i+1},\dots,U_p}(U_i) - f^i_{U_{i+1},\dots,U_p}(V_i) | \leq C \tr_N\left((U_i-V_i)(U_i-V_i)^*\right)^{1/2}.$$
	Thus thanks to Corollary 4.4.28 from \cite{alice} we have that,
	\begin{align*}
	\P\left( \left| f(U) - \E_W\left[f(W U)\right]  \right|\geq \delta \right) &\leq \sum_i \P\left( \left| f^i_{U_{i+1},\dots,U_p}(U_i) - \E_{Y_i}\left[ f^i_{U_{i+1},\dots,U_p}(W_i U_i) \right] \right| \geq \frac{\delta}{p} \right) \\
	&\leq 2p\ e^{- \left(\frac{\delta}{p C}\right)^2 N} .
	\end{align*}
	
	\noindent Finally we conclude by taking the real and imaginary part of $f$.
	
\end{proof}

We can now prove the concentration inequality that we will use in the rest of this paper. To simplify notations we will write $M$ instead of $M_N$. We also set $ Z^{NM}=(Z^N\otimes I_M, I_N\otimes Y^M) $ and $Z=(z\otimes 1, 1\otimes y)$.

\begin{prop}
	\label{2lips}
	
	Let $P\in \PP_d$, there are polynomials $H_P, K_P\in\R^+[X]$ which only depends on $P$ such that for any $N,M$,
	\begin{align*}
	\P\Bigg( &\left|\ \norm{\widetilde{P}(U^N\otimes I_M, Z^{NM})} - \E\left[\norm{\widetilde{P}(U^N\otimes I_M, Z^{NM})}\right]\ \right| \geq \delta + \frac{K_P(\norm{Z^{NM}})}{N} \Bigg) \leq e^{-\frac{\delta^2 N}{H_P\left(\norm{Z^{NM}}\right)}} ,
	\end{align*}
	
	\noindent where $\norm{Z^{NM}} = \sup\limits_i \norm{Z_i^{NM}}$.
	
\end{prop}

\begin{proof}
	
	\noindent We set $G_N : X \mapsto \norm{\widetilde{P}(X\otimes I_M, Z^{NM})}$. One can find a polynomial $L_P\in\R^+[X]$ such that for any $N$ and $Z^{NM}$, $$ |G_N(X) - G_N(Y)| \leq L_P\left( \norm{Z^{NM}} \right) \sum_i \norm{X_i - Y_i } , $$
	
	\noindent where $\norm{.}$ is the operator norm. Besides
	$$ \sum_i \norm{X_i - Y_i} \leq \sum_i \tr_N\left( (X_i - Y_i)^* (X_i - Y_i) \right)^{1/2} . $$
	
	\noindent Hence with Proposition \ref{2concentr}, there is a polynomial $H_P\in\R^+[X]$ which only depends on $P$ such that for any $N,M$,
	\begin{align*}
	\P\Big( &\left|\ \norm{\widetilde{P}(U^N\otimes I_M, Z^{NM})} - \E_W\left[\norm{\widetilde{P}(W U^N\otimes I_M, Z^{NM})}\right]\ \right| \geq \delta \Big) \leq e^{-\frac{\delta^2 N}{H_P\left(\norm{Z^{NM}}\right)}} .
	\end{align*}
	
	\noindent Besides for any matrix $U\in \U_N$, there exist $S\in SU_N$ and $\theta\in [0,2\pi]$ such that $U = e^{\i \frac{\theta}{N}} S$. Indeed we just have to pick $\theta$ such that $e^{\i \theta} = \det(U)$. Thus there is a polynomial $K_P$ such that 
	\begin{equation*}
	\left| \E_W\left[\norm{\widetilde{P}(W U^N\otimes I_M, Z^{NM})}\right] - \E\left[\norm{\widetilde{P}(U^N\otimes I_M, Z^{NM})}\right] \right| \leq \frac{K_P(\norm{Z^{NM}})}{N} .
	\end{equation*}
	
	\noindent This concludes the proof.
	
\end{proof}

We can now prove Theorem \ref{2strongconv}. Firstly, we can assume that $Z^{N}$ and $Y^M$ are deterministic matrices by Fubini's theorem. The convergence in distribution is a well-known theorem, we refer to \cite{alice}, Theorem 5.4.10. We set $g$ a function of class $\mathcal{C}^{\infty}$ which takes value $0$ on $(-\infty,1/2]$ and value $1$ on $[1,\infty)$, and belongs  to $[0,1]$ otherwise. Let us define $f_{\varepsilon}:t\mapsto g\left(\varepsilon^{-1} \left(t - \norm{\widetilde{P}\widetilde{P}^*(u\otimes 1, Z)}\right)\right)$. By Theorem \ref{2imp1}, there exists a constant $C$ which only depends on $P$, $\sup_M \norm{Y^M}$ and $\sup_N \norm{Z^{N}}$ (which is finite thanks to the strong convergence assumption on $Y^M$ and $Z^{N}$) such that,

\begin{align*}
\Bigg| &\E\left[\tr_{MN}\Big(f_{\varepsilon}\left(\widetilde{P}\widetilde{P}^*\left(U^N\otimes I_M,Z^{NM}\right)\right)\Big)\right] - MN \tau_N\otimes\tau_M\Big(f_{\varepsilon}\left(\widetilde{P}\widetilde{P}^*\left(u\otimes I_M,Z^{NM}\right)\right)\Big) \Bigg| \\
&\leq C \varepsilon^{-6} \frac{M^3\ln^2(N)}{N} .
\end{align*}

According to Theorem A.1 from \cite{male}, $(u,Z^{N})_{N\geq 1}$ converges strongly in distribution towards $(u,z)$ since, given a system of free semi-circular variable, we can write $u_i = f(x_i)$ for a specific function $f$ built with the help of Lemma 3.1 of \cite{collins_male}. Besides thanks to Lemma \ref{2tensorconv} and Corollary 17.10 from \cite{exact}, we have that $(u\otimes I_{M},1\otimes Y^{M})_{M\geq 1}$ converges strongly in distribution towards $(u\otimes 1, 1\otimes y)$. In Theorem \ref{2strongconv}, we are interested in the situation where $Z^{NM} = Z^N\otimes I_M$ or $Z^{NM} = I_N\otimes Y^M$. So, without loss of generality, we restrict ourselves to this kind of $Z^{NM}$. We know that $(u\otimes I_M, Z^{NM})$ converges strongly towards $(u\otimes 1, Z)$, but since the support of $f_{\varepsilon}$ is disjoint from the spectrum of $\widetilde{P}\widetilde{P}^*(u\otimes 1,Z)$, thanks to Proposition \ref{2hausdorff}, for $N$ large enough, $\tau_N\otimes\tau_M\Big(f_{\varepsilon}\left(\widetilde{P}\widetilde{P}^*\left(u\otimes I_M,Z^{NM}\right)\right)\Big) = 0$ and therefore,
\begin{equation}
\label{2Neps}
\E\left[\tr_{MN}\Big(f_{\varepsilon}\left(\widetilde{P}\widetilde{P}^*\left(U^N\otimes I_M,Z^{NM}\right)\right)\Big)\right] \leq C \varepsilon^{-6} \frac{M^3\ln^2(N)}{N} .
\end{equation}

\noindent  Hence, we deduce for $N$ large enough,

\begin{align*}
&\E\left[\norm{\widetilde{P}\widetilde{P}^*\left(U^N\otimes I_M,Z^{NM}\right)}\right] - \norm{\widetilde{P}\widetilde{P}^*(u\otimes 1,Z)} \\
&\leq \varepsilon + \int_{\varepsilon}^{\infty}\P\left( \norm{\widetilde{P}\widetilde{P}^*\left(U^N\otimes I_M,Z^{NM}\right)} \geq \norm{\widetilde{P}\widetilde{P}^*(u\otimes 1,Z)} + \alpha \right)\  d\alpha \\
&\leq \varepsilon + \int_{\varepsilon}^{K} \P\left( \tr_{NM}\left(f_{\alpha}\left(\widetilde{P}\widetilde{P}^*\left(U^N\otimes I_M,Z^{NM}\right)\right)\right) \geq 1 \right)\  d\alpha \\
&\leq \varepsilon + C'\varepsilon^{-6} \frac{M^3\ln^2(N)}{N} .
\end{align*}

\noindent Finally we get that,
$$ \limsup_{N\to \infty}\ \E\left[\norm{\widetilde{P}\widetilde{P}^*\left(U^N\otimes I_M,Z^{NM}\right)}\right] \leq \norm{\widetilde{P}\widetilde{P}^*(u\otimes 1,Z)} .$$

\noindent Thanks to Proposition \ref{2lips}, by taking $\delta_N = N^{-1/4}$ and using Borel-Cantelli lemma, we get that almost surely,
$$ \lim_{N\to \infty} \norm{\widetilde{P}\widetilde{P}^*\left(U^N\otimes I_M, Z^{NM}\right)} - \E\left[\norm{\widetilde{P}\widetilde{P}^*\left(U^N\otimes I_M,Z^{NM}\right)}\right] = 0 $$

\noindent And consequently almost surely,
$$ \limsup_{N\to \infty} \norm{\widetilde{P}\widetilde{P}^*\left(U^N\otimes I_M, Z^{NM}\right)} \leq \norm{\widetilde{P}\widetilde{P}^*(u\otimes 1,Z)} . $$

Besides, we know thanks to Theorem 5.4.10 of \cite{alice} that if $h$ is a continuous function taking positive values on $\left(\norm{\widetilde{P}\widetilde{P}^*(u\otimes 1, Z)}-\varepsilon, \infty \right)$ and taking value $0$ elsewhere. Then 
$$\frac{1}{MN}\tr_{MN}(h(\widetilde{P}\widetilde{P}^*(U^N\otimes I_M, Z^{NM})))$$
converges almost surely towards $\tau_{\A}\otimes_{\min}\tau_{\B} (h(\widetilde{P}\widetilde{P}^*(u\otimes 1, Z)))$. If this quantity is positive, then almost surely for $N$ large enough so is $\frac{1}{MN}\tr_{MN}(h(\widetilde{P}\widetilde{P}^*(U^N\otimes I_M, Z^{NM})))$, thus

$$ \norm{\widetilde{P}\widetilde{P}^*(U^N\otimes I_M, Z^{NM})} \geq \norm{\widetilde{P}\widetilde{P}^*(u\otimes 1, Z)} - \varepsilon .$$

\noindent Since $h$ is non-negative and the intersection of the support of $h$ with the spectrum of $\widetilde{P} \widetilde{P}^*(u\otimes 1, Z)$ is non-empty, we have that $h(\widetilde{P}\widetilde{P}^*(u\otimes 1, Z)) \geq 0$ and is not $0$. Besides, we know that the trace on the space where $z$ is defined is faithful, and so is the trace on the $\mathcal{C}^*$-algebra generated by a single semicircular variable, hence by Theorem \ref{2freesum}, so is $\tau_{\A}$. Thus, since both $\tau_{\A}$ and $\tau_{\B}$ are faithful, by Lemma \ref{2faith}, so is $\tau_{\A}\otimes_{\min}\tau_{\B}$ and $\tau_{\A}\otimes_{\min}\tau_{\B} (h(\widetilde{P}\widetilde{P}^*(u\otimes 1, Z)))>0$. As a consequence,  almost surely,
$$ \liminf_{N\to \infty} \norm{\widetilde{P}\widetilde{P}^*\left(U^N\otimes I_M, Z^{NM}\right)} \geq \norm{\widetilde{P}\widetilde{P}^*(u\otimes 1,Z)} . $$

\noindent We finally conclude  thanks to the fact that for any $y$ in a $\CC^*$-algebra, $\norm{yy^*} = \norm{y}^2$.

\section*{Acknowledgements}

The author would like to thanks his PhD supervisors Beno\^it Collins and Alice Guionnet for proofreading this paper and their continuous help, as well as Mikael de la Salle for its advices and the proof of Lemma \ref{2normineq}. The author was partially supported by a MEXT JASSO fellowship and Labex Milyon (ANR-10-LABX-0070) of Universit\'e de Lyon.

\vspace{1cm}

\bibliographystyle{abbrv}


\end{document}